%% file: Loewner.tex
\documentclass[12pt]{article}
\usepackage{amssymb}

\usepackage{amsmath,amsthm,bbm}
\usepackage{version}
\usepackage{color}

\input def_latex

\newcommand\Ha{\mathbb{H}}

\newcommand\Pid{{{\Pi}^d}}

\renewcommand\L{{\mathcal L}}
\renewcommand\P{\mathcal P}
\renewcommand\S{\mathcal S}

\newcommand\diag{\mathrm {diag~}}
\renewcommand\phi{\varphi}
\renewcommand\epsilon{\varepsilon}

\def\bpt{$B$-point \ }
\def\bpn{$B$-point}
\def\bpp{$B$-point. \ }
\def\bps{$B$-points}
\def\cpt{$C$-point \ }

\def\cps{$C$-points}
\renewcommand\ln{\lambda_n}
\renewcommand\DD{{\mathcal D}}
\newcommand\aalpha{a}
\newcommand\lnde{\L^d_n(E)}
\newcommand\lne{\L_n(E)}
\newcommand\one{\mathbbm{1}}
\newcommand\mn{M_n}
\newcommand\smn{SAM_n^d}
\newcommand\csmn{CSAM_n^d}

\newcommand\csa{CSA^d}
\newcommand\lbe{\L_\partial(E)}
\renewcommand\Re{\mathrm{Re\, }}
\renewcommand\Im{\mathrm{Im\, }}
\newcommand\nt{\stackrel{\rm nt}{\to}}

\numberwithin{equation}{section}
\title{Operator Monotone Functions and L\"owner Functions of Several Variables}
\author{Jim Agler
\thanks{Partially supported by National Science Foundation Grant
DMS 0801259}
\and
John E. M\raise.5ex\hbox{c}Carthy
\thanks{Partially supported by National Science Foundation Grant DMS 0966845}
\and
N. J. Young
\thanks{Partially supported by London Mathematical Society Grant 4918}}
\date{June 30, 2013}

\begin{document}

\bibliographystyle{plain}

\maketitle
\abstract{
We prove generalizations of L\"owner's results on matrix monotone
functions to several variables. We give a characterization of when a function
of $d$ variables is locally monotone on $d$-tuples of commuting self-adjoint $n$-by-$n$
matrices. We prove a generalization to several variables of Nevanlinna's theorem 
describing analytic functions that map the upper half-plane to itself and satisfy a growth condition.
We use this to characterize all rational functions of two variables that are operator monotone.

}
\input introduction

\input notation

\input model

\input omitvalue

\input pick
\input nevrep

\input matmonSep13
\input locopmon

\input glob


\bibliography{Loewner.bbl}

\end{document}

%% file: def_latex.tex
\usepackage{latexsym}
\usepackage{amssymb}
\usepackage{euscript}
\usepackage[dvips]{graphics}
\usepackage{epsf}

\let\cal=\mathcal      

\def\mcc{M\raise.5ex\hbox{c}C}
\def\mccarthy{M\raise.5ex\hbox{c}Carthy}

\def\eg{{\it e.g. }}
\def\ie{{\it i.e. }}

\def\h{{\cal H}}

\def\M{{\cal M}}
\def\N{{\cal N}}

\def\ga{\gamma}
\def\De{\Delta}
\def\La{\Lambda}

\def\a{\alpha}
\def\l{\lambda}

\def\vare{\varepsilon}

\let\i=\infty

\def\la{\langle}
\def\ra{\rangle}
\def\={\ = \ }

\def\A{{\cal A}}

\def\DD{D_\delta(k,w)}

\def\C{\mathbb C}
\def\R{\mathbb R}
\def\T{\mathbb T}
\def\D{\mathbb D}

\def\inn{\ \in \ }

\def\dis{\displaystyle}

\def\be{\setcounter{equation}{\value{theorem}} \begin{equation}}
\def\ee{\end{equation} \addtocounter{theorem}{1}}
\def\beq{\begin{eqnarray*}}
\def\eeq{\end{eqnarray*}}
\def\se{\setcounter{equation}{\value{theorem}}} 
\def\att{\addtocounter{theorem}{1}}
\def\vs{\vskip 5pt}
\def\bs{\vskip 12pt}

\def\bp{{\sc Proof: }}
\def\ep{{}{\hfill $\Box$} \vskip 5pt \par}

\def\bl{\begin{lemma}}
\def\el{\end{lemma}}
\def\bt{\begin{theorem}}
\def\et{\end{theorem}}
\def\bprop{\begin{prop}}
\def\eprop{\end{prop}}
\def\bd{\begin{definition}}
\def\ed{\end{definition}}
\def\br{\begin{remark}}
\def\er{\end{remark}}
\def\bexer{\begin{exercise}}
\def\eexer{\end{exercise}}
\def\bfig{\begin{figure}}
\def\efig{\end{figure}}

\newtheorem{theorem}{Theorem}[section]
\newtheorem{prop}[theorem]{Proposition}
\newtheorem{lemma}[theorem]{Lemma}
\newtheorem{cor}[theorem]{Corollary}

\newtheorem{definition}[theorem]{Definition}

\theoremstyle{definition}
\newtheorem{example}[theorem]{Example}
\newtheorem{remark}[theorem]{Remark}

%% file: introduction.tex
\section{Introduction}
\label{secint}

In 1934, K. L\"owner published a very influential paper \cite{lo34}
studying functions on an open interval $E \subseteq \R$
that are matrix monotone,
\ie functions $f$ with the property that whenever $S$ and $T$ 
are self-adjoint matrices whose spectra are in $E$ 
then
\be
\label{eqa1}
S \leq T \qquad \Rightarrow \qquad f(S) \leq f(T) .
\ee
This property is equivalent (see Subsection~\ref{subseciglo}) to being 
locally matrix monotone, \ie if $S(t)$ is a $C^1$ arc of self-adjoint
matrices with $\sigma(S(t)) \subset E$ 
then
\be
\label{eqa1.1}
S'(t) \geq 0 \qquad \Rightarrow \qquad 
\dis \frac{d}{dt} f(S(t))  \geq 0 .
\ee

Roughly speaking, L\"owner showed that if one fixes a dimension $n$ and
wants (\ref{eqa1}) or (\ref{eqa1.1}) to hold for $n$-by-$n$ self-adjoint matrices,
then  certain  
matrices derived from the values of $f$
must all be positive semi-definite. 
As $n$ increases, the conditions become more restrictive.
In the limit as $n \to \i$ 
(equivalently, if one passes to
self-adjoint operators on an infinite dimensional Hilbert space), then
a necessary and sufficient condition is that
the function $f$ 
must have an analytic continuation to a function
$F$ that maps the upper half-plane 
to itself.

\vs
The goal of this paper is to extend the above notions to several variables. In particular, we want to study functions of $d$ variables applied
to $d$-tuples of commuting self-adjoint operators. Given two $d$-tuples
$S= (S^1,\dots, S^d)$ and $T=(T^1, \dots, T^d)$, we shall say that
$S \leq T$ if and only if $S^r \leq T^r$ for every $1 \leq r \leq d$.
We want to study functions that satisfy (\ref{eqa1}) or (\ref{eqa1.1})
for $d$-tuples.

%
%
%
\vs
Before we can describe our results, we must first give a more detailed description
of the one-dimensional case. We recommend the book \cite{don74} by
W. Donoghue for a well-written account from a modern perspective.
See also the paper \cite{roro75}.

Note that there is another approach to extending L\"owner's results to several
variables where the operators $S^1, \dots, S^d$ act on
different spaces $\h^1, \dots, \h^d$, and
$f(S)$ is interpreted to act on $\h^1 \otimes \dots \otimes \h^d$.
We refer the reader to the papers \cite{han03,siva02,kor61} 
and references therein.


\subsection{Dimension one}
\label{subseci1}

Let $E$ be an open set in $\R$, and let $n \geq 2$ be a natural number.
The L\"owner class $\L_n^1(E)$ is the set of $C^1$ functions $f : E \to \R$
with the property that, whenever $\{ x_1, \dots, x_n \}$
is a set of $n$ distinct points in $E$,
then the matrix $A$,
defined by
\[
A_{ij} \ = \
\left \{ \begin{array}{lcl} \displaystyle 
\frac{ f(x_j) - f(x_i)}{x_j - x_i }
 & \mbox{ if } & i \neq j \\
f'(x_i)
 & \mbox{ if } & i = j ,
  \end{array} \right.
   \]
is positive semi-definite.

We shall let $M_n$ denote the $n$-by-$n$ complex matrices, $SAM_n$
the self-adjoint $n$-by-$n$ matrices, and $SA$ the bounded self-adjoint
operators on an infinite dimensional separable Hilbert space.

\bd
\label{defi9}
A function $f$ 
is locally $n$-matrix monotone on the open set $E \subset \R$ if,
whenever
$S$ is
in $SAM_n$
with $\sigma(S)$ consisting of $n$ distinct points
in  $E$, and
$S(t)$ is a $C^1$ curve
in $SAM_n$
with $S(0) = S$ and
$\dis \frac{d}{dt} S(t) |_{t=0} \geq 0 $,
then
$\dis \frac{d}{dt} f(S(t)) |_{t=0} \geq 0 $.
\ed

\begin{remark}
\label{rema1}
This definition is slightly different from the one in the first paragraph, where the eigenvalues
were not required to be distinct. We use this definition to be consistent with
the multivariable Definition~\ref{defi10} below. However,
using formula (6.6.31) in \cite{horjoh91} for 
$\dis \frac{d}{dt} f(S(t))$, it is easy to show that in the one variable
case the two different definitions are equivalent\footnote{
This formula says that for a $C^1$ arc $S(t) = U(t) \Lambda(t) U^*(t)$, with $U(t)$ unitary and
$\Lambda(t)$ diagonal with diagonal entries $\lambda_1(t), \dots, \lambda_n(t)$, and a $C^1$
 function $f$, then 
\[\frac{d}{dt} f(S(t)) \ =\
U(t) \left( [ \Delta f(\lambda_i(t), \lambda_j(t)) ] \circ [ U(t)^* S'(t) U(t) ] \right) U(t)^* ,\]
where $\Delta$ means the matrix of divided differences, and $\circ$ denotes the Schur product.
}.
\end{remark}


We shall say that $f$
is $n$-matrix monotone on $E$, or $M_n$-monotone,
if, whenever $S$ and $T$
are in $SAM_n$
and all their eigenvalues lie in $E$, then 
(\ref{eqa1}) holds.
To emphasize the difference from locally monotone, we shall
also call $n$-matrix monotone functions {\em globally $M_n$-monotone}.
Replacing $SAM_n$ by $SA$, we get the definitions
of locally operator monotone and operator monotone.

\bt 
[L\"owner]
Let $E \subseteq \R$ be open, and let $f \in C^1(E)$. Then $f$ is 
locally $n$-matrix monotone on $E$
if and only if $f$ is in $\L_n^1(E)$.
\label{thmi1}
\et
We shall use $\Pi$ to denote the upper half-plane, $\{ z \in \C : \Im z > 0 \}$.

\bd
Let $E \subseteq \R$ be open.
The Pick class on $E$, denoted $\P(E)$, is the set of real-valued functions
$f$ on $E$
for which there exists 
an analytic function $F : \Pi \to \overline{\Pi}$ such that
$F$ extends analytically across $E$ and\footnote{
The notation $ y \searrow 0$ means $y$ decreases to $0$.
The notation $r \nearrow 1$ means $r$ increases to $1$.}
$$
\lim_{y \searrow 0} F(x +iy) = f(x) \ \ \forall \ x \inn E .
$$
\ed

\bt 
[L\"owner]
Let $E \subseteq \R$ be open, and let $f \in C^1(E)$. The following are equivalent:

(i) The function  $f$ is locally operator monotone on $E$.

(ii) 
The function 
$f$ is in $\L_n^1(E)$ for all $n$.

(iii) 
The function 
$f$ is in $\P(E)$.

\label{thmi2}
\et

\subsection{Dimension $d \geq 2$ : Local results}

We shall let $\csmn$ denote the set of $d$-tuples of commuting 
self-adjoint $n$-by-$n$ matrices, and $\csa$ be the set of
$d$-tuples of commuting self-adjoint bounded operators. If $S$ is 
a commuting $d$-tuple of self-adjoint operators
acting on the Hilbert space $\h$,
and $f$ is a real-valued continuous (indeed, measurable) function on the joint spectrum
of $S$ in $\R^d$, then $f(S)$ is a well-defined self-adjoint operator
on $\h$.
 
\bd
\label{defilocop}
Let $E$ be an open set in $\R^d$, and $f$ be a real-valued
$C^1$ function on
$E$.
Say $f$ is locally operator monotone on $E$ if, whenever
$S$ is
in $\csa$
with $\sigma(S) \subset E$,
and
$S(t)$ is a $C^1$ curve
in $\csa$
with $S(0) = S$ and
$\dis \frac{d}{dt} S(t) |_{t=0} \geq 0 $,
then
$\dis \frac{d}{dt} f(S(t)) |_{t=0}$ exists and is $ \geq 0 $.
\ed
We shall not concern ourselves in this paper
on what conditions on $f$ guarantee that $f(S(t))$ is differentiable;
for these see \eg \cite{pe85}.


\bd
\label{defi10}
Let $E$ be an open set in $\R^d$, and $f$ be a real-valued
$C^1$ function on
$E$.
We say $f$ is locally $\mn$-monotone on $E$ if, whenever
$S$ is
in $\csmn$
with $\sigma(S) = \{ x_1,\dots, x_n \}$
consisting of $n$ distinct points in  $E$, and
$S(t)$ is a $C^1$ curve
in $\csmn$
with $S(0) = S$ and
$\dis \frac{d}{dt} S(t) |_{t=0} \geq 0 $,
then
$\dis \frac{d}{dt} f(S(t)) |_{t=0}$ exists and is $ \geq 0 $.
\ed

We define the L\"owner classes in $d$ variables, $\lnde$,
by:
\bd
Let $E$ be an open subset of $\R^d$. The set
$\lnde$ consists of all real-valued $C^1$-functions on $E$
that have the following property:
whenever $\{x_1,\dots,x_n\}$ are $n$ distinct points in $E$,
there exist positive semi-definite $n$-by-$n$ matrices
$A^1,\dots,A^d$
so that
\beq
A^r(i,i) &\=&
 \left.
 \frac{\partial f}{\partial x^r} \right|_{x_i}
\\
 {\rm and\ }\quad f(x_j) - {f(x_i)} &\=&
 \sum_{r=1}^d
 (x^r_j -  x^r_i) A^r(i,j) \qquad \forall\
 1\leq i,j \leq n .
 \eeq
 \ed

Here is our  $d$-variable version of 
Theorem~\ref{thmi1}; we prove it as Theorem~\ref{thmma}.

\bt

Let $E$ be an open set in $\R^d$, and $f$ a real-valued $C^1$ function
 on $E$.
 Then $f$ is locally $\mn$-monotone if and only if $f$ is in $\lnde$.
\et

\vs
In generalizations of Theorem~\ref{thmi2}, there turns out to be
a difference between the case $d=2$ and $d>2$.

\bd
The L\"owner class, $\L^d$, is the set of 
functions 
$F: \Pid \to \overline \Pi$ 
with the property that there exist $d$ positive 
semi-definite kernel
functions $A^r, \ 1 \leq r \leq d,$
on $\Pid$ such that
\[
F(z) - \overline{F(w)} \= (z^1 - \bar w^1) A^1 (z,w) 
\ + \ \dots \ +  \ (z^d - \bar w^d) A^d (z,w) .
\]
\ed           
When $d=1$ or $2$, the L\"owner class coincides with the set of all
analytic functions from $\Pid$ to $\overline \Pi$, but for $d \geq 3$
it is a proper subset (see Remark~\ref{rem216}). 

\bd
\label{defi11}
Let $E \subseteq \R^d$ be open.
The class $\L(E)$ is the set of real-valued functions
$f$ on $E$
for which there exists
an analytic function $F$ in $\L^d$ such that 
$F$ extends analytically across $E$ and
$$
\lim_{y \searrow 0} F(x^1 +iy, \dots, x^d +iy) = f(x^1, \dots, x^d) \ \ \forall \ x \inn E .
$$
\ed

We prove the following result as Theorem~\ref{thmloa}.

\bt
Let $E$ be an open set in $\R^d$, and $f$ a 
real-valued $C^1$ function on $E$.
The following are equivalent:

(i) The function $f$ is locally operator monotone on $E$.

(ii) The function $f$ is in $\lnde$ for all $n$.

(iii) The function $f$ is in $\L(E)$.
\et

\subsection{Local to Global}
\label{subseciglo}

In one variable, provided $E$ is an interval, 
 local monotonicity implies global
monotonicity immediately.
Indeed, suppose $S \leq T$, and let $S(t) = (1-t) S + t T $.
Then $S'(t) = T-S \geq 0$, so
\be
\label{eqi4}
f(T) - f(S) 
\=
\int_{0}^1 
\frac{d}{dt} f(S(t))\  dt \ \geq 0.
\ee

If $E$ is not convex, this argument fails. Indeed, the function $-1/x$ is locally $n$-matrix monotone
on $\R \setminus \{ 0 \}$ for all $n$; but it is only globally monotone on sets that lie entirely
on one side of $0$. A result of Chandler \cite{cha76} says that functions that are globally operator 
monotone on a set $E$ always extend to be globally monotone on the convex hull of $E$.


For intervals, (\ref{eqi4}) shows 
that the word ``locally'' can be dropped in both Theorem~\ref{thmi1} and \ref{thmi2}. 
One problem in going to several variables is that this simple argument 
no longer works, because 
one may not be able to connect $S$ and $T$ by a 
path of commuting $d$-tuples.
Indeed, the following example shows that there need not be 
any commuting tuples between two given ones.

\begin{example}
\label{exi1}
Let $S$ and $T$ be pairs in $CSAM^2_2$ given by
\beq
S &=&
\left( 
\left(
\begin{array}{cc}
0&0 \\
0&5 \end{array} \right) 
\, , \,
\left(
\begin{array}{cc}
1&0 \\
0&0 \end{array} \right) \right) \\
T &=&
\left( 
\left(
\begin{array}{cc}
4&2 \\
2&6 \end{array} \right) 
\, , \,
\left(
\begin{array}{cc}
2&2 \\
2&4 \end{array} \right) \right) .
\eeq
If $R$ is in $CSAM^2_2$ and $S \leq R \leq T$,
it can be shown that either $R=S$ or $R=T$.
\end{example}

%
%

We have been unable to 
resolve the question of whether
the $n$-matrix monotone functions on a connected open set $E$ are a 
proper subset of the locally $n$-matrix monotone functions on $E$.
However, as $n$ tends to infinity and we pass to locally operator monotone
functions, analyticity enters the picture, and makes the problem more tractable ---
see Subsection~\ref{subsecie}.

\subsection{The Nevanlinna Representation}
To prove $(iii) \Rightarrow (i)$ in Theorem~\ref{thmi2}, one must understand
analytic functions that map the upper half-plane to itself. A key fact
is a characterization due to R.~Nevanlinna \cite{nev22}
which says that, provided they have some regularity
at infinity, they are all Cauchy transforms
of measures on the line.
\bt
[Nevanlinna] If $F: \Pi \to \Pi$ is analytic and satisfies
\[
\limsup_{y \to \i} \,  y \, | F(iy) - C | < \i ,
\]
for some $C \in \R$, 
then there exists a unique finite positive Borel measure $\nu$ on $\R$ so that
\be
\label{eqi14}
F(z) \=
C +  \int \frac{d\nu(t)}{t-z} .
\ee
\label{thmi3}
\et
Nevanlinna's theorem was used by M.~Stone to prove the spectral theorem
\cite{sto32},
but one can adopt the reverse viewpoint, and write (\ref{eqi14}) in terms
of the resolvent of a self-adjoint. Indeed, let $X$ be the self-adjoint operator
of multiplication by the independent variable on $L^2(\nu)$, and 
$v$ the vector in $L^2(\nu)$ that is $1$ a.e. 
Then (\ref{eqi14}) can be rewritten as 
\be
\label{eqi15}
F(z) \=
C +  \la  (X-z)^{-1} v, v \ra.
\ee
This representation turns out to be useful in 
studying operator monotonicity,
because then 
\be
F(S) \=
C I +  R_v^*  (I\otimes X-S \otimes I)^{-1}
 R_v,
\label{eqi16}
\ee
where $R_v: \mathcal{H} \to \mathcal{H \otimes M}$ is given
by $R_v: \xi\mapsto \xi\otimes v$.

There is a several variable analogue of Theorem~\ref{thmi3}.
It may require first perturbing $F$.

\bd
\label{defi7}
For each real number $t$, define
\[
\rho_t(z) = \frac{z+t}{1-tz} .
\]
For $F \inn \L^d$, define $F_t$ by
\[
F_t(z^1,\dots,z^d)  \ :=\ \rho_{t} \circ F ( \rho_t (z^1),\dots, \rho_t(z^d)) .
\]
\ed

The following theorem follows from Theorem~\ref{thmd4}.
We shall say that a function $F$ on $\Pid$ is analytic on a neighborhood
of infinity if the function $F(1/z^1,\dots, 1/z^d)$ extends to be analytic on a neighborhood
of the origin.
In Theorem \ref{thmd4}, a weaker assumption is placed on $F$ than being analytic
in a neighborhood of infinity. 
\bt
\label{thmi5}
Let $F$ be in $\L^d$, and 
assume that $F$ is analytic in a neighborhood of infinity.
Then for all sufficiently small $t$, with at most countably
many exceptions, the function $F_t$ has the following representation.
There is a Hilbert space $\M$, a densely defined 
self-adjoint operator $X$ on $\M$,
a vector $v$ in $\M$, a real constant $C$,
and $d$ orthogonal projections $P^1, \dots, P^d$ with $\sum_{r=1}^d P^r = I_\M$
so that
\be
\label{eqi17}
F_t(z) \=
C +  \la  (X-\sum_{r=1}^d z^r P^r )^{-1} v, v \ra.
\ee
\et

\subsection{Dimension $d \geq 2$ : Global operator monotonicity}
\label{subsecie}

If $E$ is an open set in $\R^d$, we shall say that a real-valued function $f$ defined on $E$
is {\em globally operator monotone}, or just operator monotone for short, if, for every $n$, whenever
$S$ and $T$ are in 
$\csmn$, with $S \leq T$, and the joint eigenvalues of both $S$ and $T$ lie in $E$, then
\[
f(S) \leq f(T) .
\]

Using 
the representation (\ref{eqi17}), we can prove results on (global) operator
monotonicity.
With notation as in Theorem~\ref{thmi5}, let us say that the
$\mu$-resolvent of $X$ is the set of points $$
\{ (z^1, \dots, z^d) \inn \C^d
\, : \,
X-\sum_{r=1}^d z^r P^r \ \ {\rm has\ a\ bounded\ inverse} \}.
$$

We prove the following result as Theorem \ref{thmg1}.

\bt
Let $X$ be a densely-defined self-adjoint 
operator on a Hilbert space
$\M$, let $v$ be a vector in $\M$, let $C$ be a real constant,
and let $P^1, \dots, P^d$ be projections with orthogonal ranges that
sum to the identity.
Let $F$ be given by
\[
F(z) \= 
C +  \la  (X-\sum_{r=1}^d z^r P^r )^{-1} v, v \ra.
\]
Let $E$ be an open box in $\R^d$ 
that is in the $\mu$-resolvent of $X$.
Then $F$ is globally operator monotone on $E$.
\et

As an application, we can give a complete characterization of
the rational functions of two variables that are operator monotone on 
rectangles. This is Theorem \ref{thmg2}.

\bt
Let $F$ be a rational function of two variables. 
Let $\Gamma$ be the zero-set of the denominator
of $F$. Assume $F$ is real-valued on $\R^2 \setminus \Gamma$.
Let $E$ be an open rectangle in $\R^2 \setminus \Gamma$. Then $F$ is globally operator
monotone on $E$ if and only if $F$ is in $\L(E)$, that is if and only if $F$ is the restriction to
$E$ of an analytic function from $\Pi^2$ to $\overline{\Pi}$ that extends analytically across $E$.
\et

%% file: notation.tex
\section{Some Notation}
\label{secnot}

We shall let $\D$ denote the unit disk in the complex plane, $\Pi$ the upper half-plane
$\{ z \, : \, \Im(z) > 0 \}$, and $\Ha$ the right half-plane
$\{ z \, : \, \Re(z) > 0 \}$.
We shall let
\be
\label{rev68}
\a(\l) \= i \frac{1+\l}{1-\l}
\ee
be a linear fractional map that maps $\D$ to $\Pi$, and
\be
\label{rev69}
\beta(z) \= \frac{z-i}{z+i}
\ee
be its inverse.

We shall let $d$ denote the number of variables.
If $z$ is a point in $\Pi^d$, we shall use $z^1,\dots,z^d$ to
denote its components; likewise $\l = (\l^1,\dots,\l^d)$ will be a point in
$\D^d$. We shall write $S= (S^1,\dots,S^d)$ for a $d$-tuple of matrices or
operators, and use $\| S \|$ for $\max_{1\leq r \leq d} \| S^r \|$.
We shall also use $\alpha$ and $\beta$ to denote the maps from $\D^d$ to $\Pid$ and back again
that are defined by applying $\alpha$ and $\beta$ coordinate-wise.

A kernel on a set $E$ is a map $K: E \times E \to \C$ 
with the property that for every finite set $\{\l_1,\dots,\l_N \}$ of distinct points
in $E$, the matrix $[K(\l_j,\l_i)]$ is positive semi-definite.

\bd
The Pick class, $\P^d$, is the set of analytic functions $F : \Pid \to \overline{\Pi}$.
\ed
\bd
The Schur class, $\S^d$, is the set of analytic functions $\phi: \D^d \to \overline{\D}
.$
\ed
\bd
The Carath\'eodory class, ${\mathcal C}^d$, is the set of analytic functions $\psi : \D^d \to \overline{\Ha}$.
\ed
\bd
The L\"owner class, $\L^d$, is the set of functions $F: \Pid \to \overline \Pi$
with the property that there exist $d$ kernel functions $A^r, \ 1 \leq r \leq d$
on $\Pid$ such that
\be
\label{eqn51}
F(z) - \overline{F(w)} \= (z^1 - \bar w^1) A^1 (z,w) \ + \ \dots \ + 
\ (z^d - \bar w^d) A^d (z,w) .
\ee
\ed
\bd
The Schur-Agler class, $\A^d$, is the set of functions $\phi: \D^d \to \overline \D$
with the property that there exist $d$ kernel functions $B^r, \ 1 \leq r \leq d$
on $\D^d$ such that
\be
\label{eqn52}
1 - \phi(\l) \overline{\phi(\mu)}
\=
(1 - \l^1 \bar \mu^1) B^1(\l,\mu) \ +
\ \dots \ + \ 
(1 - \l^d \bar \mu^d) B^d(\l,\mu) .
\ee
\ed

When the dimension is clear, we shall drop the superscript $d$.
\begin{remark}
If we
exclude the constant function $1$ from $\S$,  we have the identification
\be
\label{eqn1}
F \inn \P \quad \iff \quad \beta \circ F \circ \a \inn \S
\quad \iff \quad -i F \circ \alpha \inn {\mathcal C}.
\ee
Moreover, we also have (again excluding the constant function $1$)
\be
\label{eqn53}
F \inn \L \quad \iff \quad \beta \circ F \circ \a \inn \A,
\ee
(see Lemma~\ref{propp1}).
As all our results are trivial for constant functions, we shall 
use (\ref{eqn1}) and
(\ref{eqn53}) without explicitly mentioning
the exclusion of the constant function $1$.
\end{remark}
The following change of variables formula is in \cite{kv04}. A function is in $\A^d$ if and only if it
is analytic and maps 
$d$-tuples of commuting strict contractions to contractions; a function is in $\L^d$ if and only if it
is analytic and maps
$d$-tuples of commuting operators 
with strictly positive imaginary parts\footnote{
We say an operator $T$ has strictly positive imaginary part if there exists $ \alpha > 0$
such that $(T - T^*)/2i \geq \alpha I$.} 
to operators with positive imaginary parts.

\begin{lemma}
\label{propp1}
The function $F: \Pid \to \C$ is in the L\"owner class
if and only if $\phi := \beta \circ F \circ \alpha$ is in the Schur-Agler class
 $\A^d$.
\end{lemma}
\bp
Define $\phi = \beta \circ F \circ \a.$
Then $\phi$ is in $\A^d$ if and only if there are kernels $B^r$ on $\D$ such that
\be
\label{eqp2}
1 - \phi(\l) \overline{\phi(\mu)}
\=
\sum_{r=1}^d
(1 - \l^r \bar \mu^r) B^r(\l,\mu) .
\ee
When $z = \a(\l)$ and $w = \a(\mu)$, (\ref{eqp2}) becomes
\be
\label{eqp3}
1 - \beta\circ F (z) \overline{\beta \circ F (w)} \=
\sum_{r=1}^d \left( 1 - \left[ \frac{z^r -i}{z^r +i} \right]
\overline{
\left[ \frac{w^r - i}{w^r +i}\right] } \right)
B^r(\beta(z), \beta(w)) .
\ee
Rearranging (\ref{eqp3}), we get
\be
\label{eqp4}
F(z) - \overline{F(w)}
\=
\sum_{r=1}^d
(z^r - \bar w^r) \,
\frac{F(z) +i}{z^r +i}
\frac{\overline{F(w)} -i}{\overline{w^r} -i}
B^r(\beta(z),\beta(w)) .
\ee
If $A^r$ is defined for $r=1,\dots,d$ by
$$
A^r(z,w)
\=
\frac{F(z) +i}{z^r +i}
\frac{\overline{F(w)} -i}{\overline{w^r} -i}
B^r(\beta(z),\beta(w))
$$
(\ref{eqp4}) becomes
\be
\label{eqp1}
F(z) - \overline{F(w)} \=
\sum_{r=1}^d (z^r - \bar w^r) A^r (z,w) ,
\ee
which means $F$ is in $\L^d$. Reversing the argument gives the converse.
\ep

\begin{remark}
\label{rem216}
It is known that $\A^d = \S^d$ for $d = 1$ or $2$, and that for $d \geq 3$
$\A^d \subsetneq \S^d$ \cite{cradav,var74,ag90}. It follows similarly that the 
L\"owner class equals the Pick class in dimensions
$1$ and $2$, and is strictly contained in it for $d \geq 3$.
By 
Theorem 5.5.1 of \cite{rud69}, rational inner functions are dense in the unit ball
of 
$\S^d$ in the topology of uniform convergence on compacta. Therefore there must
be rational inner functions in $\S^d \setminus \A^d$ for $d \geq 3$. 
By (\ref{eqn1}) and (\ref{eqn53}), it follows  
that for each $d \geq 3$,
there are rational 
functions that are real on $\R^d$ and that are in $\P^d \setminus \L^d$.
\end{remark}

%% file: model.tex
\section{Models, \bps\ and \cps}
\label{secmod}

For a function $\phi$ in $\A^d$, we can take the representation (\ref{eqn52})
and decompose the $B^r$'s as Gramians to get a Hilbert space model for $\phi$.
That means we find
a separable Hilbert space $\mathcal{M}$, an orthogonal decomposition
 of $\mathcal{M}$,
\be
\label{eqmo1}
\mathcal{M} = \mathcal{M}^1 \oplus \cdots \oplus \mathcal{M}^d,
\ee
and an analytic map $u: \mathbb{D}^d \rightarrow \mathcal{M}$ such that
\be
\label{1.3}
1 - \overline{\phi(\mu)}\phi(\lambda) = \sum_{r=1}^d \ (1 - \overline{\mu^r} \lambda^r)
\langle  u_\lambda^r,u_\mu^r
\rangle_{\M^r}
\ee
 for all $\lambda, \mu \in \mathbb{D}^d$, where
 we write
$u_\lambda$ for $u(\lambda)$, $P^r$ for the
projection onto $\M^r$, and
$u_\lambda^r $ for $ P^r [u_\lambda ]$.

We shall view (\ref{eqmo1}) interchangeably
as a graded Hilbert space (\ie one with a given orthogonal decomposition) 
or as a single Hilbert space with $d$ given projections $P^1, \dots, P^d$ that
are orthogonal and
add up to the identity.

In general, if $\eta \in \mathcal{M}$, we set $\eta^r = P^r [\eta]$.
If $\lambda \in \mathbb{C}^d$,
we may regard $\lambda$ as an operator on $\mathcal{M}$ by letting
\be
\label{eqp13}
\lambda \eta = \lambda^1 \eta^1 + \dots + \lambda^d \eta^d.
\ee
Equation (\ref{1.3}) can then be rewritten as
\be
\label{eql8}
1 - \overline{\phi(\mu)}\phi(\lambda) = \langle  (1 - \mu^* \lambda) u_\lambda, u_\mu \rangle .
\ee

A lurking isometry argument yields the following result
\cite{ag90}.
\begin{theorem} \label{thml3}
If $(\mathcal{M}, u)$ is a model of $\phi \in \A^d$, then there
exist  $a\in\C$, vectors $\beta, \gamma \in \mathcal{M}$ and a
linear operator $D: \mathcal{M} \rightarrow \mathcal{M}$ such that the operator
\beq\label{abcd}
\left[ \begin{array}{cc} a & 1\otimes\beta \\ \ga\otimes 1 & D \end{array} \right]
\eeq
is a contraction on $ \C \oplus \mathcal{M}$ and,  for all $ \l\in\D^d$,
\se\att
\begin{eqnarray}
\label{eqaa6}
 (1-D \lambda) u_\lambda &=& \gamma,\\
   \phi(\l)  &=& a + \langle \l u_\l, \beta  \rangle .
\label{eqaa7}
\end{eqnarray}
 \end{theorem}
\att
With notation as in Theorem \ref{thml3}, we shall call $(a,\beta,\gamma,D)$ a realization
of $(\mathcal{M},u)$.

One can rewrite (\ref{eqaa6}) and (\ref{eqaa7}) as 
\be
\label{rev21}
\phi(\lambda) \= a + \beta^* \lambda(I - D \lambda)^{-1} \gamma .
\ee
How one can go from (\ref{rev21}) to (\ref{1.3}) is discussed in
\cite{baltre98} and \cite{babo10}.

If we start instead with the representation (\ref{eqn51}) of a function $F$ in $\L^d$,
we can decompose the kernels $A^r$ as the Grammians of some vectors $v^r$,
in auxiliary separable 
Hilbert spaces $\N^r$. Then we get, in the analogous notation to above,
\begin{eqnarray}
\nonumber
F(z) - \overline{F(w)} &\=&
\sum_{r=1}^d (z^r - \bar w^r) A^r (z,w)
\\
\nonumber
&=& \sum_{r=1}^d (z^r - \bar w^r) \la v^r_z , v^r_w \ra_{\N^r} \\
\label{eqaa9}
&=&  \la (z - w^*)v_z , v_w \ra_{\N}.
\end{eqnarray}
\att
This decomposition leads to a lurking self-adjoint argument, which
we shall discuss in Section~\ref{secnev}.

\bs

In \cite{amy10a}, we introduced the concept of a \bpt for $\S$.
Let us give a unified definition for each of the classes $\S, \P$ and $\mathcal C$;
notice that it depends on the codomain of the function.
\bd Let $U$ and $V$ be fixed domains, and $f: U \to \overline{V}$
an analytic function. A point $\tau$ in $\partial U$ is called a
\bpt of $f$ if there is a sequence $\l_n$ of points in $U$ that converge
to $\tau$ and such that
\be
\label{eqn2}
\frac{{\rm dist}(f(\ln), \partial V)}
{{\rm dist}(\ln, \partial U)}
\ee
is bounded.
\ed

So, for example,
a point $\tau$ in $\partial \Pid$ is a \bpt for a function $F$ in $\P^d$
if there exists some sequence $z_n$
in $\Pid$ that tends to $\tau$ and such that the quantity
$$
\frac{\Im F(z_n)}{ \min_{r\in\{1,\dots,d\}} (\Im z^r_n )}
$$
is bounded.

For a function in $\L^d$ (respectively, $\A^d$) we shall call a point $\tau
$ a \bpt
if it is a \bpt for the function thought of as an element of $\P^d$ (resp.
$\S^d$).

For each of the three classes $\S, \P,$ and $\mathcal C$,
it follows from results of F. Jafari \cite{jaf93}
and M.~Abate \cite{ab98}
that if
$\tau$ is a \bpn, then
the ratio (\ref{eqn2}) remains bounded for every sequence $\l_n$
that tends to $\tau$ non-tangentially. Moreover, the function
$f$ will then have a non-tangential limit at $\tau$.
(A sequence $\l_n$ in $U$ tends to the point $\tau$ non-tangentially
if $\l_n$ tends to $\tau$ and
$$
\frac{{\rm dist}(\ln, \tau )}
{{\rm dist}(\ln, \partial U)}
$$
is bounded.)

The following result
was proved in \cite{amy10a} for $d=2$, but the proof  generalizes to any $d$.
We shall need it in the proof of Theorem~\ref{thmd3}.
\begin{lemma}
\label{leml3}
Let $\phi \in \A^d$ and $\tau \in \T^d$. Let $(\mathcal{M}, u)$ be a model
for $\phi$, and  $(a,\beta,\gamma,D)$ be a realization. The
following are equivalent.

(i) $\tau$ is a \bpt for $\phi$.

(ii) For some  sequence $\ln$ converging to $\tau$ non-tangentially,
the sequence $\| u_{\ln} \|$ is bounded.

(iii) For any sequence $\ln$ converging to $\tau$ non-tangentially,
the sequence $\| u_{\ln} \|$ is bounded.

(iv) The vector $\gamma$ is in the range
of $(I - D \tau)$.

Moreover, if $u_{\ln}$ converges to a vector weakly as $\ln$ tends to $\tau$
non-tangentially, then $u_{\ln}$ converges in norm.
The vector $u_\tau := \lim_{r \nearrow 1} u_{r \tau}$ exists for every
\bpt $\tau$.
\end{lemma}

A stronger condition than being a \bpt is being
a $C$-point. 
\bd
\label{defmod5}
A point $x \in \R^d$ is a $C$-point
for $F \inn \L$ if 
there are complex numbers $\eta^1, \dots ,\eta^d$
and a real number $c$ so that
\be
\label{eqref6}
F(z) - c - \sum_{r=1}^d \eta^r (z^r - x^r)
\=
o(\| z - x \|)
\ee
as $z$ tends to $x$ non-tangentially.
\ed
In particular, if $F$ is differentiable at $x$ and $F(x)$ is real,
then $x$ is a \cpt for $F$, as (\ref{eqref6}) then holds as $z$ tends to $x$
from any direction.

The following result was proved in \cite{amy10a}.
\bprop
\label{propmodc}
Suppose $F \inn \L$ has a model as in (\ref{eqaa9}).
If $x$ is a \cpt for $F$, then as $z$ converges to $x$  non-tangentially
from $\Pid$, the vectors $v_z$ converge in norm, to some
vector $v_x$ in $\N$.
\eprop

%% file: omitvalue.tex
\section{Analytically continuing Pick functions}
\label{secomit}

Suppose $F$ is analytic on $\Pid$, and 
$E$ is an open set in $\R^d$. What conditions on $F$ guarantee
that it can be analytically continued across $E$?
The edge-of-the-wedge theorem (see Theorem~\ref{thmc5} below)
is a common tool to give such extensions. 
Checking the hypotheses, however, 
requires knowledge of the values of $F$ as one approaches points of $E$ not
just non-tangentially but also tangentially.
If $F$ is in the Pick class $\P^d$, Theorem~\ref{thmc3}
says that it suffices to know that every point
of $E$ is  a \bpt (which can be checked by looking at the
values of $F$ on the inward-pointing normal).

As we are using bars to denote closure, we shall use stars 
for the complex conjugate of a set,
and write $\Pi^*$ for the lower half-plane.

\subsection{One dimension}
\label{subsecc1}

To understand the situtation, let us first consider the one dimensional case.
Let $\psi : \D \to {\Ha}$ be non-constant. Then $\psi$ has a Herglotz representation;
if we assume $\psi(0)$ is positive, then
\be
\label{eqc1}
\psi(z) \= \int_{0}^{2\pi} \frac{e^{i\theta} + z}{e^{i\theta} -z} d \mu(\theta)
\ee
for some positive measure $\mu$.
There is an elegant analysis of when $\psi$ has \bps\ in the paper 
\cite{sar98c} by D.~Sarason,
where the following two propositions are proved.
Proposition~\ref{propc1} is originally due to M.~Riesz
\cite{rie31}, and Proposition \ref{propc2} to
 R.~Nevanlinna \cite{nev29b}.
\bprop
\label{propc1}
Let $\psi$ be given by (\ref{eqc1}),
and let $\tau$ be a point in $\T$. Then $\psi$ has a \bpt
at $\tau$ 
if and only if 
\be
\label{eqc2}
\int \frac{1}{| e^{i\theta} - \tau |^2} d \mu(\theta)
\ < \ \i .
\ee
\eprop
If $\phi =\beta \circ ( i \psi)$ is the Cayley transform of $\psi$,
there is a distinction between \bps\ where $\phi(\tau)$ 
equals $1$, corresponding to $\psi(\tau) = \i$, and all other cases.
\bprop
\label{propc2}
Let $\dis \phi = \frac{\psi -1}{\psi + 1}$, where $\psi$ is
given by (\ref{eqc1}),
and let $\tau$ be a point in $\T$.
Then $\phi$ has a \bpt at $\tau$ with $\phi(\tau) \neq 1$ if 
and only if (\ref{eqc2}) holds.
The function $\phi$ has a \bpt at $\tau$ with $\phi(\tau) = 1$ if 
and only if $\tau$ is a mass point of $\mu$.
\eprop

Suppose now that $\phi : \D \to \D$ has an open arc $I$ of \bps. Can $\phi$ be extended
analytically across $I$? If we know that $\phi$ omits a value on $I$, then the answer is yes.
Indeed, after a M\"obius map,  we can assume that $\phi$ is the Cayley transform of some $\psi$ as (\ref{eqc1}).
If condition (\ref{eqc2}) holds on an open arc $I$, then $\mu$ must vanish on $I$ by Lemma~\ref{lemc1} below.
But then the formula (\ref{eqc1}) gives an analytic function on the extended plane less $\T \setminus I$.

However, without the assumption that $\phi$ omits a value, the answer may be no, as 
Example~\ref{exc1} 
shows.

\begin{lemma}
\label{lemc1}
Suppose $\mu$ is a measure on $[-\pi,\pi)$ and (\ref{eqc2}) holds for $\tau = e^{ix}$ for every $x$  in an open arc $I \subset [-\pi,\pi)$.
Then $\mu(I) = 0$.
\end{lemma}
\bp
For $\mu$ a.e. point $x$ in $I$, there is a constant $ c > 0$ such that
$$
\mu [ x - \frac{1}{k}, x + \frac{1}{k}]
\ \geq \ c \frac{2}{k} ,
$$
by 
\cite[Thm 8.6 and Thm 8.10]{rudrc}.
For such an $x$, we have
\beq
\int_{-\pi}^{\pi} \frac{1}{| e^{i\theta} - \tau |^2} d \mu(\theta) 
&\
\geq \ &
\int_{x - 1/k}^{x+1/k} \frac{1}{| e^{i\theta} - \tau |^2} d \mu(\theta) \\
&\ \geq \ &
\frac{1}{| 1 - e^{i/k} |^2} \, \mu[
 x - \frac{1}{k}, x + \frac{1}{k}]
 \\
& \geq & c k .
\eeq
Letting $k$ tend to infinity, the integral would be infinite; so $\mu$ 
must put no mass on $I$.
\ep
\begin{example}
\label{exc1}
Here is an example of a function in the Schur class of one variable
that has \bps\ at every 
point of $\T$ but that cannot be analytically continued across every arc.

Let $\tau_n = e^{ix_n}$ be a sequence in $\T$ that converges to $1$.
Let $\dis c_n = 2^{-n} | 1 - \tau_n |^2$.
Then for every $\l \in \T$, the quantity
$\dis \left| \frac{1-\tau_n}{\l - \tau_n} \right| $ is less than or equal to $1$
 for all but finitely many values of $n$.
Therefore
\be
\label{eqc3}
\sum 2^{-n}
\left| \frac{1-\tau_n}{\l - \tau_n} \right|^2 \ < \ \i
\ee
for every $\l$.

Let $\dis \mu = \sum 2^{-n} \delta_{x_n}$, let $F$ be the Herglotz transform of $\mu$, and let  
$\dis \phi = \frac{F -1}{F + 1}$. By Proposition \ref{propc2}, we have that $\phi$ has every point of $\T$ as a \hbox{\bpn,}
but $\phi$ cannot be analytically continued across any arc containing $1$, as it takes the value $1$ infinitely often
on any such arc.
\end{example}

\subsection{$d$ dimensions}

Our goal is to prove the following analytic continuation theorem:

\bt
\label{thmc3}
Let $E$ be an open subset of $\R^d$. 
Then there is an open set $U$ in $\C^d$ that contains $\Pid \cup E \cup {\Pi}^{*d}$ 
with the following property:
whenever  $F$ is in the Pick class, and every point of $E$ is
a \bpt for $F$, 
then there is an analytic
function $G$ on $U$ that agrees with $F$ on $\Pid$.
\et

This theorem immediately implies the omit-a-value theorem.
Let us say that a subset $E'$ of $\T^d$ is a {\em $B$-set} for $\phi$ in
${\mathcal S}^d$ if every point of $E'$ is a $B$-point for $\phi$.
\bt
\label{thmc4}
Let $E'$ be an open subset of $\T^d$.
Then there is an open set $U$ in $\C^d$ containing
$\D^d \cup E' \cup \{ \C \setminus \bar \D \}^{d}$ 
such that the following two statements are equivalent for any $\phi$
 in the Schur class:
\begin{enumerate}
\item [\rm (1)]
there is an analytic
function $\psi$ on $U$ that agrees with $\phi$ on $\D^d$;

\item[\rm (2)] the set $E'$ is a $B$-set for $\phi$ and
for every point $\tau$ in $E'$ there exists a neighborhood $V$ of $\tau$ in
$\T^d$ and a point $\omega$ in $\T$ such that no nontangential limit of
$\phi$ at any point of $V$ is equal to $\omega$.
\end{enumerate}
\et
Condition (2) says that  every point of $E'$ has a neighborhood
where  the non-tangential limits of $\phi$ omit some value in $\T$.
\bs
We need a version of the edge-of-the-wedge theorem
(Theorem C from \cite{rudeow}). We write $\R_+$ for the interval $(0,\i)$.

\bt
\label{thmc5}
[Edge-of-the-wedge]
Let $E$ be an open subset of $\R^d$.
Then there is an open set $U$ in $\C^d$ that
contains $\Pid \cup E \cup {\Pi}^{*d}$
and is such that
whenever $H$ is an analytic function on $\Pid$
with the property that
for every $g$ in
$C^\i_c (E)$,
\be
\label{eqc4}
  \lim_{\R_+^d \ni y \to 0}
  \int_E g(x) \Im H(x+iy) dx \= 0 ,
  \ee
then there is an analytic
  function $G$ on $U$ that agrees with $H$ on $\Pid$.
  \et


\bprop
\label{propc4}
Let $E$ be an open subset of $\R^d$.
Then there is an open set $U$ in $\C^d$ that contains $\Pid \cup E \cup {\Pi}^{*d}$ 
with the following property:
if $J$ is a non-empty interval in $\R$, 
$F$ is in the Pick class,  $F$ has non-tangential limits at almost every point of $E$, and these limits are all in $\R \setminus J$,
then there is an analytic
function $G$ on $U$ that agrees with $F$ on $\Pid$.
\eprop

\bp
Precomposing $F$ with a M\"obius transformation of $\Pi$ if necessary, we can assume that $J$ is an interval about
infinity, so the non-tangential limits are in some compact set $[-M,M]$ a.e.

Let $H(z) = \log ( 1 + M + F(z) )$. Then $H$ maps $\Pid$ into $\dis \{ z \in \C : 0 < \Im z < \pi \}$,
and
$$
 \lim_{y \to 0} H(x + iy) \ \in \ [ 0, \log (2M+1) ] \quad {\rm a.e.\ } x \in E .
$$

As $H$ has bounded imaginary part, we can pass the limit inside the integral
on the left-hand side of (\ref{eqc4}), and as $H$ has real boundary values,
we get that the limit is $0$.
Therefore by Theorem~\ref{thmc5} we get an analytic extension of $H$, and hence $F$,
to the desired
open set $U$.
\ep

{\sc Proof of Theorem~\ref{thmc3}:}
We can extend $F$ to $\Pi^{*d}$ by letting $F(z) = \overline {F(\bar z)}$ on $\Pi^{*d}$.
The difficulty is in showing that the definitions of $F$ on the two disjoint domains
$\Pid$ and $\Pi^{*d}$ are analytic continuations of each other across $E$.
This is a local property. If we can show that every point of $E$ has a neighborhood
on which the boundary values of $F$ take values in a bounded set, we can apply Proposition~\ref{propc4}
to conclude that the reflection of $F$ is an analytic continuation of $F$ across this 
neighborhood in $E$. Since this is true at every point, the conclusion of the theorem will follow.

For convenience, we will change variables and consider the function $\phi(\l) = (-i) F \circ \alpha$,
which is in ${\mathcal C}^d$. 

We can normalize to assume that $\phi(0,\dots,0) = 1$ and that the point of 
interest 
for $\phi$
is $\beta(0,\dots,0) = (-1,\dots,-1)$. So for some 
$ 0 < c < \frac{\pi}{5}$, the set
\be
\label{eqmay11a}
\{ ( e^{i \theta_1} , \dots, e^{i \theta_d} ) \ : \
\forall\ 1 \leq r \leq d, \  | \theta_r  | \geq \pi - 5c \}
\ee
consists of
\bps\ for $\phi$. In what follows, we shall choose $\arg$ to take values
in $[-\pi,\pi)$.

For each 
$\tau$ in the set
$$
\{ \tau \in \T^{d-1} \ : \ | \arg(\tau^j)| < 2c, \ \forall\ 1 \leq j \leq d-1 \} ,
$$
define  $g_\tau$ in $\mathcal{C}^1$ by
$$
g_\tau(z) \= \phi(z,\tau^1 z,\tau^2 z, \dots, \tau^{d-1} z ) .
$$
Then for each $\tau$, the set $$
I_{3c} \= 
\{ \sigma \in \T \, : \,  |\arg (\sigma) |  > \pi - 3c \} $$
is a set of \bps\ for $g_\tau$, and $g_\tau(0) = 1$.
Each $g_\tau$ has a Herglotz representation, and by the results of subsection~\ref{subsecc1}
the corresponding measure is supported off the set $I_{3c}$.
So
$$
g_\tau(z) \=
\int_{-\pi + 3c}^{\pi - 3c} \frac{e^{i\theta} + z}{e^{i\theta} -z} d \mu_\tau(\theta)
$$
for some probability measure $\mu_\tau$.
Therefore if $\sigma$ is in the arc $I_{c} = \{ \sigma \in \T \, : \,  
|\arg (\sigma)|  >  \pi - c \} $,
\beq
| g_\tau (\sigma) |
&\ \leq \ &
\int_{-\pi + 3c}^{\pi - 3c} \left| \frac{e^{i\theta} + \sigma}{e^{i\theta} -\sigma} \right|d \mu_\tau(\theta)
\\
&\leq&
 \sec c . 
\eeq
Therefore on the set $(I_{c})^d$ we conclude that the non-tangential limits of $\phi$ take values in
the bounded set $[ - \sec c,  \sec c]$.

Notice that $c$ does not depend on $F$: we have shown that for any $F$,
normalized to have $F(i,\dots,i) = i$, if $F$ has \bps\ on
the set $\alpha((I_{5c})^d)$, then $F$ is bounded on
$\alpha((I_{c})^d)$. By Proposition~\ref{propc4},
this latter set now has a neighborhood to which $F$ can be analytically 
extended, 
and this neighborhood can be chosen independently of $F$.
So every point $x$ in $E$ has a neighborhood $U_x$ to which all 
functions $F$ in the Pick class with \bps\  on $E$
can be extended; let $U$ be the union of all the $U_x$ as $x$ ranges over
$E$.
\ep

%% file: pick.tex
\section{The L\"owner classes}
\label{secpick}

We shall single out functions that have a representation
on subsets of $\R^d$ as in (\ref{eqn51}).

\bd
\label{deb1}
Let $E \subseteq \R^d$ be a non-empty open set,
and let $n$ be a positive integer.
We define $\lne$
to be
the set of real valued differentiable
functions
that have the following property:
whenever $\{x_1,\dots,x_n\}$ are $n$ distinct points in $E$,
there exist positive semi-definite $n$-by-$n$ matrices
$A^1,\dots,A^d$
so that
\se\att\begin{eqnarray}
\label{defb11}
A^r(i,i) &\=&
\left.
\frac{\partial f}{\partial x^r} \right|_{x_i}
\\
\nonumber
{\rm and} &&\\
\att
\label{defb12}
\quad f(x_j) - {f(x_i)} &\=&
\sum_{r=1}^d
(x^r_j -  x^r_i) A^r(i,j) \ \forall\
1\leq i,j \leq n .
\end{eqnarray}
\ed

We shall 
give an alternative description of $\L(E)$. We shall temporarily
call it $\lbe$, but we shall show in Proposition~\ref{propba1}
that it coincides with the set $\L(E)$ from Definition~\ref{defi11}.

\bd
\label{deb2}
Let $E \subseteq \R^d$ be a non-empty open set.
We shall let $\lbe$ denote
the set of differentiable real valued functions $f$
on $E$ for which
there exist positive semi-definite functions $A^1,\dots,A^d: E \times
E \to \C$ so that
\se\att\begin{eqnarray}
\label{defb14}
A^r(z,z)
&\=&
\left.  \frac{\partial f}{\partial x^r}  \right|_{z}
\\
\att
\label{defb13}
{\rm and\ }\quad
f(z) - {f(w)} &\=&
\sum_{r=1}^d
(z^r -  w^r) A^r(z,w) .
\end{eqnarray}
\ed

If $E \subseteq \R^d$, a function $f$ in
$\lbe$
can be extended to a function $F$ in $\L$ that
has $f$ as its non-tangential (which we abbreviate nt)
boundary values on $E$.

\bprop
\label{propthmp1}
Let $E \subseteq \R^d$ be open, and let $f \in \lbe$. Then there
exists $F \in \L^d$
such that every point of $E$ is a \bpt for $F$, and such that
\be
\label{eqba5}
\lim_{z \nt t} F(z) \= f(t) \qquad \forall \ t \in E .
\ee
\eprop

Proposition~\ref{propthmp1}
follows immediately from the corresponding result on the polydisk,
Theorem~\ref{thmp2}, which was proved by J.A.~Ball and V. Bolotnikov \cite{babo02}
(we are changing their language slightly; they did not explicitly use the notion
of $B$-point).

\bt
[Ball-Bolotnikov]
Let $E' \subseteq {\T^d}$ and let $\psi : E' \to \C$.
Suppose there are
positive semi-definite functions $B^1,\dots ,B^d : E' \times E' \to \C$
such that, for all $\l,\, \mu$ in $E'$,
\be
\label{eql7}
1 - \psi(\l) \overline{\psi(\mu)}
\=
(1 - \l^1 \bar \mu^1) B^1(\l,\mu) \ +
\ \dots \ + \
(1 - \l^d \bar \mu^d) B^d(\l,\mu) .
\ee
Then there is a function
$\phi$ in $\A$
such that every point of $E'$ is a \bpt for $\phi$ and
$$
 \lim_{\l \nt \tau} \phi(\l) \= \psi(\tau) \qquad \forall \ \tau \in E'
.
$$
Moreover, if $\phi$ is defined to equal $\psi$ on $E'$, the kernels $B^r$
can be extended to $E' \cup \D^d$ so that,
for all $\l,\, \mu$ in $E' \cup \D^d$,
\[
1 - \phi(\l) \overline{\phi(\mu)}
\=
(1 - \l^1 \bar \mu^1) B^1(\l,\mu) \ +
\ \dots \ + \
(1 - \l^d \bar \mu^d) B^d(\l,\mu) .
\]
\label{thmp2}
\et


We can pass back and forth between regarding functions in
$\L(E)$ as functions in the L\"owner class $\L^d$ that 
have \bps\ on $E$ (and so can be analytically extended across $E$),
and as functions that are characterized by their values on $E$
by (\ref{deb2}) 
and can then be analytically extended into $\Pi^d$.

\bprop
\label{propba1}
Let $E \subseteq \R^d$ be a non-empty open set. The following four 
sets coincide.

(i) $ \cap_{n=1}^\i \lne$.

(ii) The set $ \lbe $ defined by Definition~\ref{deb2}.

(iii) The set $\L(E)$ defined by Definition~\ref{defi11}.

(iv) The functions $f$ on $E$ for which there exists
a function $F$ in $\L^d$ such that every point $x$ of $E$
is a \bpt of $F$ and the non-tangential limit of $F$ at $x$ is $f(x)$.

\eprop
\bp
It is immediate that $(ii) \subseteq (i)$. 
Theorem~\ref{thmc3} asserts that $(iii) = (iv)$.
Proposition \ref{propthmp1} says that $(ii) \subseteq (iv)$.

To show $(iii) \subseteq (ii)$, 
choose a model for $F$ so that 
(\ref{eqaa9}) holds on $\Pid \times \Pid$:
$$
F(Z) - \overline{F(W)} \= \la (Z- W^*) v_Z, v_W \ra \qquad \forall\ Z,W \inn \Pid.
$$
As every point in $E$ is a \cpt for $F$, we can let $Z$ and
$W$ tend to points in $E$ non-tangentially, $z$ and $w$ respectively.
By Proposition~\ref{propmodc}, the vectors $v_Z$ and $v_W$
converge to $v_z$ and $v_w$. Let $$
A^r(z,w) \= \la v^r_z, v^r_w \ra
$$
and one gets (\ref{defb13}). To get (\ref{defb14}),
let $z$ be in $E$ and let $W$ in $\Pid$ tend to $z$ non-tangentially.
As $F$ is analytic at $z$, we have
\be
\label{eqwk7}
F(W) - F(z) \= \sum (W^r - z^r) 
\left. \frac{\partial f}{\partial x^r} \right|_{z}
\ + \ o(\|z-W\|) .
\ee
From the model,
\begin{eqnarray}
\nonumber
F(W) - F(z) &\=& \la (W-z) v_W, v_z \ra \\
&=& \la (W-z) v_z, v_z \ra \ + \
\la (W-z) (v_W - v_z) , v_z \ra
\label{eqwk8}
\end{eqnarray}
\att
The second term on the right of (\ref{eqwk8})
is $o(\|z - W\|)$,  
so comparing with (\ref{eqwk7}) we conclude that
$$
\left. \frac{\partial f}{\partial x^r} \right|_{z}
\= \la v^r_z, v^r_z \ra,
$$
and hence (\ref{defb14}) holds.

\vs
To prove $(i) \subseteq (ii)$,
we need 
to show that if (\ref{defb14}) and 
(\ref{defb13}) hold on every finite set,
with perhaps a different choice of $A^r$'s each time, then we can make
one choice for the 
$A^r$'s that works everywhere.

Let $f\in\cap_{n=1}^\infty \L_n(E)$.  
Consider any finite set $\{z_1,\dots,z_n\}$ of distinct points in $E$.  By
Definition \ref{deb1} there exist kernels $A^1,\dots,A^d$ on $E$ such that
equations (\ref{defb11}), (\ref{defb12}) hold, and we have
\[
A^r(i,i)\  \leq\  \displaystyle \frac{\partial f}{\partial x^r}(z_i),
  \qquad i=1,\dots,n, \ r=1,\dots,d.
\]
Since the matrices $A^r$ are positive semi-definite we also obtain bounds on the
off-diagonal entries of all the $A^r$.  Hence the set $K$ of all
$d$-tuples $(A^1,\dots,A^d)$ for which equations (\ref{defb11}), (\ref{defb12})
hold is a compact non-empty subset of $M^d_n$.

Moreover, if $(B^1,\dots,B^d)$ is a $d$-tuple of kernels on any finite
superset $Z$ of $\{z_1,\dots,z_n\}$ for which the analogs of equations
(\ref{defb11}), (\ref{defb12}) hold, then the choice of $A^r$ to be the
principal submatrix of $B^r$ corresponding to $\{z_1,\dots,z_n\}$ gives a
$d$-tuple that belongs to $K$.  
Therefore by Kurosh's theorem \cite[p.74]{arkh} or \cite[p.30]{ampi}, there is 
a $d$-tuple $(A^1,\dots,A^d)$ of kernels on $E$ such that
equations (\ref{defb11}) and (\ref{defb12}) hold for all points $z_i,\ z_j \in
E$.
\ep

%% file: nevrep.tex
\section{The $\mu$-spectral theorem}
\label{secnev}

A function in the Pick class of one variable, \ie an analytic function from $\Pi$ to $\Pi$, has an integral 
representation which can be obtained from the Herglotz representation 
(\ref{eqc1}) of functions from $\D$ to $\Ha$ by a change
of variables \cite{her11}. 

\bt
\label{thmd1} [Herglotz]
An analytic function  $F : \Pi \to \Pi$ has a unique representation of the form
\be
\label{eqd1}
F(z) \= c  \, + \, d z \, + \, \int 
\frac{1+zt}{t-z}
d \mu(t) 
\ee
where $\Im c \geq 0$ and  $d \geq 0 $, and $\mu$ is a finite
positive Borel measure on $\R$.
Conversely any function of this form is in the Pick class of one variable.
\et
If in addition $F$ decays up the imaginary axis, one gets that $F$ is the Cauchy transform of a finite measure on $\R$. This is called the Nevanlinna representation, and was proved by R.~Nevanlinna \cite{nev22}.
\bt
\label{thmd15}
[Nevanlinna] If $F: \Pi \to \Pi$ is analytic and satisfies
\be
\label{rev78}
\limsup_{y \to \i} \, | y F(iy)| < \i ,
\ee
then there exists a unique finite positive Borel measure $\nu$ on $\R$ so that 
\[
F(z) \= 
 \int \frac{d\nu(t)}{t-z} .
\]
\et

\begin{remark}
 If one considers $\psi = -i F \circ \alpha : \D \to \Ha$,
then $d$ in Theorem~\ref{thmd1} is the mass assigned to the point $1$ in the Herglotz representation of $\psi$.
Nevanlinna's condition (\ref{rev78})
in Theorem~\ref{thmd15}
is equivalent to 
saying that $\phi = \beta \circ F \circ \alpha: \D \to \D$ has a \bpt at $1$ with $\phi(1) = -1$.
\end{remark}

One can prove the spectral theorem for a (possibly unbounded)
self-adjoint operator by showing that, if $R_z$ is the resolvent, then
for any vector $u$ the function $\la R_z u, u \ra$ is in the one variable Pick class, 
and satisfies Nevanlinna's growth condition.
Then Theorem~\ref{thmd15}
gives the scalar spectral measure. 
See \cite[Chap. V]{don74} or \cite[Chap. 32]{lax02}.
Conversely, if $X$ is the operator of multiplication by the 
independent variable on $L^2(\mu)$, and $v$ is the constant function
$1$, then (\ref{eqd1}) becomes
\be
\label{eqd1a}
F(z) \= c  \, + \, d z \, + \, \la (1+zX)(X-z)^{-1} v , v \ra .
\ee

In several variables, there is also a connection between Pick functions and self-adjoint
operators, which could be called a $\mu$-spectral theorem (Theorem~\ref{thmd2} below).

\bd
Let $\M$ be a Hilbert space, with a fixed decomposition as
$\M = \M^1 \oplus \dots \oplus \M^d$. Let $T$ be a densely defined linear operator on $\M$.
For $z = (z^1,\dots,z^d)$ in $\C^d$, define the $\mu$-resolvent of $T$ at $z$ to be
$$
(T - z)^{-1} \= (T -\,  [z^1 I_{\M^1} \oplus \dots \oplus z^d I_{\M^d}]\, )^{-1} .
$$
The $\mu$-spectrum of $T$ is the complement of the set of points in $\C^d$
for which the $\mu$-resolvent is bounded.
\ed

The expressions ``$\mu$-resolvent" and ``$\mu$-spectrum" are not standard,
but they are suggested by usage in control engineering.  The notion of
$\mu$-analysis provides an approach to robust stabilization in the
presence of ``structured uncertainty" \cite{DulPag}.  Corresponding to the
projections $P^1,\dots,P^d$ one defines the ``cost function" $\mu(X)$ by
\beq
\lefteqn{
\frac{1}{\mu(X)} = \inf \{ \|T\|: T\in B(\M),}
\\&\qquad \qquad \mbox{ each } P^r\M \mbox{
reduces } T \mbox{ and } 1-TX \mbox{ is singular } \}.
\eeq


In what follows, we shall write $z$ for
$z^1 I_{\M^1} \oplus \dots \oplus z^2 I_{\M^d}$ and
$z^*$ for
$\bar z^1 I_{\M^1} \oplus \dots \oplus \bar z^d I_{\M^d}$.
Let us recall Definition~\ref{defi7}.

{\bf Definition \ref{defi7}.}
{\em
For each real number $t$, define
\be
\label{eqd15}
\rho_t(z) = \frac{z+t}{1-tz} .
\ee
For $F \inn \L^d$, define 
\[
\label{eqd596}
F_t  \ :=\ \rho_{t} \circ F \circ \rho_t .
\]
}
Note that, similarly to the maps $\alpha$ and $\beta$, we use 
$\rho_t$ on $\C^d$ to mean the component-wise action.

\bt
\label{thmd2}
Let $F$ be in $\L^d$, and $z_0$ a point in $\Pid$.
For all except at most a countable number of real numbers $t$,  
there is a Hilbert space
$\M = \M^1 \oplus \dots \oplus \M^d$, a self-adjoint operator $X$ on $\M$, a vector $v$ in $\M$,
and a real number $c$ 
so that
\be
\label{eqd2}
F_t(z) \=  c
\, + \,
\la  z v, v \ra 
\, + \,
\la (z - z_0^*) (X-z)^{-1} (z-z_0) v, 
v \ra .
\ee
Conversely, if $z_0$ is a point in $\Pid$, $c$ is a real number,
$X$ is a densely defined self-adjoint operator
on a Hilbert space  
$\M = \M^1 \oplus  \dots \oplus \M^d$, and  $v$ is a vector in $\M$, 
then
the function of $z$ given by 
the right-hand side of
(\ref{eqd2})
is in $\L^d$.
\et

\bp
($\Rightarrow$)
Let $\phi = \beta \circ F \circ \alpha$ in $\A^d$, where $\alpha$ and $\beta$
are defined in (\ref{rev68}) and (\ref{rev69}). Choose a model for
$\phi$ so that (\ref{1.3}) holds:
\se\att
\begin{eqnarray}
\nonumber
1 - \overline{\phi(\mu)}\phi(\lambda) &\=& \sum_{r=1}^d \ (1 - \overline{\mu^r} \lambda^r)
\langle  u_\lambda^r,u_\mu^r
\rangle_{\M^r}\\
\label{eqdd1}
&=& \la (1 - \mu^* \l) u_\l, u_\mu \ra_{\M}.
\end{eqnarray}
Define a linear operator $V$ by
$$
V \, : \,
\left( \begin{array}{c}
1 \\ \l u_\l
\end{array} \right)
\ \mapsto \
\left( \begin{array}{c}
\phi(\l)  \\  u_\l
\end{array} \right) ,
$$
and extend it by linearity to finite linear combinations of
vectors of the form
$$
\left( \begin{array}{c} 1 \\ \l_i u_{\l_i} \end{array} \right)
$$
where the points $\l_i$ range over $\D^d$.

$V$ is defined on a subspace of $\C \oplus \M$, and by (\ref{eqdd1}) it is isometric
on its domain.
If the codimensions of the closures of the
domain and range of $V$ are the same, $V$ can be extended to a 
unitary $U$. 
If they are different, 
after the addition of a separable infinite dimensional summand to one of the spaces 
$\M^r$, the codimensions become equal, and
one can then extend $V$ to a unitary $U$.
So we can assume that we have a unitary $U : \C \oplus \M \to \C \oplus \M$ such that
\be
\label{eqdd6}
U \, : \,
\left( \begin{array}{c}
1 \\ \l u_\l
\end{array} \right)
\ \mapsto \
\left( \begin{array}{c}
\phi(\l)  \\  u_\l
\end{array} \right) .
\ee
Now, let $\tau$ be a point in the unit circle that is {\em not} in the point spectrum
of $U$, and let 
$$
t \= -i \frac{1-\tau}{1+\tau} .
$$
As $\C \oplus \M$ is separable, the point spectrum of 
$U$ is countable, so all but countably many
real numbers $t$ will arise in this way.

Let $$
Y \= -i(U-\tau)^{-1}(U + \tau) \ : \ (U-\tau) \eta \mapsto \ -i(U+\tau) \eta .
$$
Then $Y$ is densely defined and self-adjoint. Its domain $\DD$
is $\mathrm{ran~}(U-\tau)$.
Moreover, by definition
$$
Y \, : \,
\left( \begin{array}{c}
\phi(\l) - \tau \\ (1 - \tau \l) u_\l
\end{array} \right)
\ \mapsto \
\left( \begin{array}{c}
-\, i(\phi(\l) + \tau)  \\ -\, i(1 + \tau \l) u_\l
\end{array} \right) .
$$
Therefore
\be
Y \, : \,
\left( \begin{array}{c}
1 \\ 
\\
\displaystyle{
\frac{1 - \tau \l}{\phi(\l) - \tau} u_\l
}
\end{array} \right)
\ \mapsto \
\left( \begin{array}{c}
\displaystyle{
-\, i\ \frac{\phi(\l) + \tau}{\phi(\l) - \tau}  
}
\\ 
\\
\displaystyle{
-\, i\ \frac{1 + \tau \l}{\phi(\l) - \tau} u_\l
}
\end{array} \right) .
\label{eqdd7}
\ee
Let
$$
\tilde{v}_z \= \frac{1 - \tau \beta (z)}{\beta \circ F (z) - \tau}
u_{\beta (z)} .
$$
Then one can rewrite (\ref{eqdd7}) as
\be
Y \, : \,
\left( \begin{array}{c}
1 \\ \tilde{v}_z
\end{array} \right)
\ \mapsto \
\left( \begin{array}{c}
\rho_t \circ F(z)  \\ 
- \rho_{-t} (z)  \tilde{v}_z
\end{array} \right) .
\label{eqdd8j}
\ee
Now let
$$
v_z \= \tilde{v}_{\rho_t(z)} .
$$

Then (\ref{eqdd8j}) becomes
\be
Y \, : \,
\left( \begin{array}{c}
1 \\ v_z
\end{array} \right)
\ \mapsto \
\left( \begin{array}{c}
F_t(z)  \\ 
- z v_z
\end{array} \right) .
\label{eqdd8}
\ee
As $Y$ is self-adjoint, (\ref{eqdd8}) implies that
\be
\label{eqdd9}
F_t(z) - F_t(w)^* \=  \la (z - w^*) v_z, v_w \ra .
\ee

Let $X$ be the compression of $-Y$ to $\M$.
By Lemma~\ref{lemd1}, $X$ is self-adjoint with
dense domain equal to $\DD \cap \M$.

If $\ga$ is in $\DD \cap \M$, then 
$$
Y \,
\left( \begin{array}{c} 0 \\ \ga \end{array} \right)
\=
\left( \begin{array}{c} L (\ga)  \\ -X \gamma \end{array} \right) .
$$
for some linear functional $L$.

Define $ v = v_{z_0}$, and let $\aalpha = F_t(z_0)$. Then
$$
Y \,
\left( \begin{array}{c} 1 \\  v \end{array} \right)
\=
\left( \begin{array}{c} \aalpha \\ - z_0 v \end{array} \right) .
$$
For $z \in \Pid$ let
\beq
Y \,
\left( \begin{array}{c} 1 \\  v_z \end{array} \right)
&=&
Y \,
\left( \begin{array}{c} 1 \\  v \end{array} \right)
\ +\
Y \,
\left( \begin{array}{c} 0 \\ v_z - v \end{array} \right)
\\
&=& 
\left( \begin{array}{c} \aalpha \\ - z_0 v \end{array} \right)
\ + \
\left( \begin{array}{c} L (v_z -v)  \\ -X (v_z -v) \end{array} \right) .
\eeq
By (\ref{eqdd8}), we get the equations
\begin{eqnarray}
\label{eqd6}
\att
F_t(z) & \=& \aalpha + L(v_z - v)
\\
\att
\label{eqd7}
zv_z &\=& z_0 v + X(v_z -v).
\end{eqnarray}
For $\ga \in \DD \cap \M$,
\beq
\la \left( \begin{array}{c} 0 \\ \ga \end{array} \right),
\left( \begin{array}{c} \aalpha \\  - z_0 v \end{array} \right) \ra
&=&
\la \left( \begin{array}{c} 0 \\ \ga \end{array} \right),
Y
\left( \begin{array}{c} 1 \\  v \end{array} \right) \ra
\\
&=& \la Y \left( \begin{array}{c} 0 \\ \ga \end{array} \right),
\left( \begin{array}{c} 1 \\  v \end{array} \right) \ra
\\
&=& \la \left( \begin{array}{c} L(\ga) \\ -X \ga \end{array} \right),
\left( \begin{array}{c} 1 \\  v \end{array} \right) \ra.
\eeq
Therefore
\be
\label{eqd8}
L(\ga) \=
- \la \ga, z_0 v \ra + \la X \ga, v \ra .
\ee
If $z$ is in $\Pid$, then by Lemma~\ref{lemd2}, $X-z$ is invertible, so
(\ref{eqd7}) yields
\be
\label{eqd9}
v_z - v \=
(X-z)^{-1} (z-z_0) v .
\ee
Combining equations (\ref{eqd6}) to (\ref{eqd9}), we get
\begin{eqnarray}
\att \se \nonumber
F_t(z) &=&
\aalpha - \la v_z - v , z_0 v \ra + \la zv_z - z_0 v , v \ra 
\\
\nonumber
&=&
\aalpha - \la v_z - v, z_0 v -z^* v \ra  +
\la (z - z_0 ) v, v \ra \\
\nonumber
&=& 
\aalpha 
+ \la (X-z)^{-1} (z-z_0) v , (z^* - z_0) v \ra
 \\
\label{eqd10}
&&
\ +\  \la (z - z_0 ) v, v \ra. 
\end{eqnarray}
By (\ref{eqdd9}),
$$
\aalpha - \bar \aalpha \= 
F_t(z_0) - \overline{F_t(z_0)}
\= 
  \la (z_0 - z_0^*) v, v \ra ,
$$
so $c := \aalpha - \la z_0 v, v \ra$ is real.
Then (\ref{eqd10}) becomes 
(\ref{eqd2}), as desired.

\vs

$(\Leftarrow)$ To prove the converse, suppose $X$ is a self-adjoint operator
on $\M$ with dense domain $\DD'$. Let  $F(z)$ be given by
the right-hand side of (\ref{eqd2}).
Define $v_z$ by (\ref{eqd9}), \ie
\be
\label{eqd11}
v_z 
\= v + (X-z)^{-1} (z-z_0) v .
\ee
Define a linear functional $L$ on $\DD'$ by
\[
L(\ga) \=
- \la \ga, z_0 v \ra + \la X \ga, v \ra .
\]
Let $\DD$ be the linear span in $\C \oplus \M$ of the vector
$ \left( \begin{array}{c} 1 \\  v \end{array} \right)$
and the vector space $ 0 \oplus \DD'$. Let $a = c + \la z_0 v,v \ra$.
Finally, define $Y$ on $\DD$ by
\[
 Y \,
\left( \begin{array}{c} t \\ \ga + tv \end{array} \right)
\=
t \,\left( \begin{array}{c} \aalpha \\ - z_0 v \end{array} \right)
\ + \
\left( \begin{array}{c} L (\ga)  \\ - X (\ga) \end{array} \right) 
\]
for $t$ in $\C$.
It is routine to verify that $Y$ is symmetric.
Moreover, by (\ref{eqd11}), $(v_z - v)$ is in the domain of $X-z$, 
and therefore in $\DD'$.
So $ \left( \begin{array}{c} 1 \\ v_z \end{array} \right)$
is in $\DD$ for every $z$ in $\Pi$, and therefore
\be
\label{eqd16}
\la Y
\left( \begin{array}{c} 1 \\ v_z \end{array} \right),
\left( \begin{array}{c} 1 \\ v_w \end{array} \right) \ra
\=
\la
\left( \begin{array}{c} 1 \\ v_z \end{array} \right),
Y \left( \begin{array}{c} 1 \\ v_w \end{array} \right) \ra.
\ee
Expanding (\ref{eqd16}) and rearranging, one gets
\[
F(z) - \overline{F(w)} \=  \la (z - w^*) u_z, u_w \ra ,
\]
and hence 
$F$ is in $\L^d$.
\ep

\begin{lemma}
\label{lemd1}
With notation as in the proof of the forward direction of Theorem~\ref{thmd2},
the domain of $X$ is $\DD \cap \M$. This domain is dense, and
the operator $X$ is self-adjoint.
\end{lemma}
\bp
Since $\DD$ is dense in $\C \oplus \M$, there are vectors $\xi_n$ in $\M$ that converge to
zero and such that $\left( \begin{array}{c} 1 \\ \xi_n \end{array} \right)$
are in $\DD$. If $\gamma$ is any vector in $\M$, there are vectors 
$\left( \begin{array}{c} a_n \\ \eta_n \end{array} \right)$
in $\DD$ that converge to $ \left( \begin{array}{c} 0 \\ \gamma \end{array} \right)$,
hence so do the vectors 
$$
\left( \begin{array}{c} a_n \\ \eta_n \end{array} \right) \ - \ 
a_n \, \left( \begin{array}{c} 1 \\ \xi_n \end{array} \right) \=
\left( \begin{array}{c} 0 \\ \eta_n - a_n \xi_n \end{array} \right).
$$
Therefore $\DD \cap \M$ is dense in $\M$.

Let $P$ be the projection from $\C \oplus \M$ onto $\M$.
Let $X =-  P Y |_\M$ with domain $\DD' = \DD \cap \M$.
Then for $\gamma,\ \eta$ in $\DD'$, we have
\beq
\la X \gamma, \eta \ra &\=& - \la P Y \ga, \eta \ra
\\
&=& - \la Y \ga, \eta \ra \\
&=& - \la \ga, Y \eta \ra \\
&=& \la \ga , X \eta \ra .
\eeq
So $X$ is symmetric.

To prove $X$ is self-adjoint, it remains to show that $X$ and $X^*$ have the same
domain. Assume that there is some vector $\eta$ in $\M$ such that
\[
| \la X \ga , \eta \ra | \ \leq \ C \| \ga \| 
\]
for all $\ga \in \DD'$.
Then for every vector $\left( \begin{array}{c} c \\ \delta \end{array} \right)$
in $\DD$ of norm at most one, we have
\beq
| \la Y \left( \begin{array}{c} c \\ \delta \end{array} \right),\,
\left( \begin{array}{c} 0 \\ \eta \end{array} \right) \ra |
&=&
|\la Y \left[ \left( \begin{array}{c} 0 \\ \delta - c \xi_1 \end{array} \right) \, 
+ \,
c\left( \begin{array}{c} 1 \\ \xi_1 \end{array} \right)\right]
,
 \left( \begin{array}{c} 0 \\ \eta \end{array} \right) \ra |
\\
&\leq&
C \| \delta - \xi_1 \| + |c| \, \| Y \left( \begin{array}{c} 1 \\ \xi_1 \end{array} \right) \|
\, \| \eta \| \ \leq \ C'.
\eeq
So $\left( \begin{array}{c} 0 \\ \eta \end{array} \right)$ is in $\DD$, and 
therefore $\eta$ is in $\DD'$.
\ep

\begin{lemma}
\label{lemd2}
Let $X$ be a densely defined self-adjoint operator on $\M$.
The $\mu$-spectrum of $X$ is disjoint from $\Pid \cup \Pid^*$.
Moreover,
$$
\| (X-z)^{-1} \| \ \leq \ 
1/ \min_{1 \leq r \leq d} ( |\Im z^r|) \qquad \forall \, z \in \Pid
\cup \Pid^* . 
$$
\end{lemma}
\bp
Let $X$ be self-adjoint on $\M = \M^1 \oplus \dots \oplus \M^d$, 
and let $z = (x^1 + i y^1, \dots, x^d + i y^d)$ be a point in $\Pid
\cup \Pid^*$.
Then for any $v = v^1 \oplus \dots \oplus v^d$ in $\M$,
\[
\la (X - z) v, v \ra 
\=
\la (X -x ) v, v \ra  -i ( y^1 \| v^1 \|^2 + \dots + y^d \| v^d \|^2 ) 
. \]
The first summand on the right is real, so $X-z$ is bounded below by
$\min( |y^r|)$, and therefore has a left inverse.
Applying the same argument to $z^*$, we get that $X-z^*$ has a left inverse,
and taking adjoints we get that $X-z$ has a right inverse also.
\ep

When $F$ decays at infinity, we can sharpen Theorem~\ref{thmd2} to get a theorem like 
Nevanlinna's Theorem~\ref{thmd15}. (It has long been known that the $\limsup$ in Nevanlinna's theorem can be replaced by a $\liminf$).

Let us write $\one$ for $(1,1,\dots,1)$, and $s\one$ for $(s,s,\dots,s)$, etc.


\bt
\label{thmd3}
Let $F$ be in $\L^d$, and assume $F$ has a representation
as in (\ref{eqd2}) with $t=0$. Then the following are equivalent.

\begin{enumerate}
\item[\rm(i)]
$$ 
\liminf_{y \to \i} \ y |F(iy \one)|  \, < \, \i .
$$

\item[\rm(ii)]
There exists a vector $v_1$ in $\M$
so that
\be
\label{eqd17}
F(z) \=  \la (X-z)^{-1} v_1 , v_1 \ra 
\qquad z \inn \Pid .
\ee


\item[\rm(iii)]
The function $\phi = \beta \circ F \circ \alpha $ in $\S$ has a \bpt at $\one$ 
and
\newline
$\phi(\one) = -1 $.


\item[\rm(iv)]
$\lim_{y \to \i} F(iy\one) = 0$ and the vectors
$v_{(iy\one)}$
defined by {\rm (\ref{eqd11})} satisfy
$$
 \liminf_{y \to \i} \  y \| v_{(iy\one)} \| \, < \, \i .
$$

\item[\rm(v)]
The vector $v$ is in the domain of $X$ and $\lim_{y \to \i} F(iy\one) = 0$.
\end{enumerate}
\et
 
\bp
$\rm(i) \Rightarrow (iii)$
Let $\phi = \beta \circ F \circ \alpha$.
Condition (i) becomes
\be
\label{eqd18}
\liminf_{s \to 1} \ \frac{1+s}{1-s} \left| \frac{1 + \phi(s\one)}{1 - \phi(s\one)} \right|
\ < \ \i .
\ee
The left-hand side of (\ref{eqd18}) dominates 
\[
\frac{1- |\phi(s\one)|}{1-s} ,
\]
so $\one$ is a \bpp In order for (\ref{eqd18}) to hold, we must
have $\phi(\one) = -1$.

$\rm(iii) \Leftrightarrow (iv)$ As the proof of Proposition~\ref{propp1} shows, one can 
pass between a 
model $({\M }, u)$ for $\phi$ and a model $(\M , v)$ for $F$ 
by letting
\se\att
\begin{eqnarray}
\label{eqd19}
v^r_z &\=& \left( \frac{F(z) + i }{z^r + i}\right) \ u^r_{\beta(z)} \\
\nonumber
u^r_\l &=& \left( \frac{1 - \phi(\l)}{1 - \l^r} \right) \  v^r_{\alpha(\l)} \qquad r=1,\dots,d.
\end{eqnarray}  

By Lemma~\ref{leml3},
$\phi$ having a \bpt at $\one$ is equivalent to  
 $u_{(r\one)}$ being bounded as $r \to 1^-$. 
Moreover, $\phi(r\one)$ tending to $-1$ 
is the same 
as $ F(iy\one)$ 
tending to $0$ as $y \to \i$.
And as long as $F(iy\one)$ has a finite limit, 
(\ref{eqd19}) says that $u_{(r\one)}$ is bounded iff
$[y \, v_{(iy\one)}]$ is.

$\rm(iv) \Rightarrow (v)$ As $X$ is densely defined and self-adjoint, it is closed.
By (\ref{eqd7}), the vectors $v_z - v$ all lie in $\DD'$, the domain of $X$.
Let $z = (iy\one)$ and let $y \to \i$. Then $v - v_{(iy\one)}$ tends to $v$.
Moreover, $X(v - v_{(iy\one)}) = z_0 v -iy v_{(iy\one)}$
contains a bounded sequence as $y \to \i$, and therefore a 
subsequence that converges weakly to some vector, $w$ say.
So $(v,w)$ is in the weak closure of the graph of $X$, therefore in the
graph of $X$, and hence $v$ is in
$\DD'$.

%

$\rm(v) \Rightarrow (ii)$
If $v \in \DD'$, then (\ref{eqd8}) becomes
\be
\label{eqd20}
L(\ga) = \la \ga , Xv -  z_0 v  \ra .
\ee
Let $v_1 = (X-z_0) v $.
Then (\ref{eqd7}) says 
\be
\label{eqd21}
(X-z) v_z = v_1.
\ee
Combining (\ref{eqd6}), (\ref{eqd20})  and (\ref{eqd21}), we get
\be
\label{eqd22}
F_t(z) \= a - \la v, v_1 \ra + \la (X-z)^{-1} v_1, v_1 \ra .
\ee
Now let $z = (iy\one)$ in (\ref{eqd22}) and let $y \to \i$. By Lemma \ref{lemd2},
the last term on the right tends to zero, so we must have $ 
a - \la v, v_1 \ra = 0$.

$\rm(ii) \Rightarrow (i)$ Lemma \ref{lemd2} implies that 
\[
\| F(z) \| \ \leq \ \| v_1 \|^2  / \min ( | \Im z^r | ),
\]
and so $\rm(i)$ follows.
\ep

%
%
%
%
%
%

For later use, let us record a slight variant of Theorem~\ref{thmd3}; it is
proved in the same way.
\bt
\label{thmd4}
Let $F$ be in $\L^d$, and assume $F$ has a representation
as in (\ref{eqd2}) with $t=0$. Then the following are equivalent.

\begin{enumerate}

\item[\rm(i)]
There exists a constant $C \in \R$ so that
$$  \qquad\qquad  \liminf_{y \to \i} \  y |F(iy \one) \, - \, C| \, < \, \i .
$$

\item[\rm(ii)]
There exists a vector $v_1$ in $\M$
and a constant $C$ in $\R$
so that
\[
F(z) \=  C \, + \, \la (X-z)^{-1} v_1 , v_1 \ra 
\qquad z \inn \Pid .
\]

\item[\rm(iii)]
The function $\phi = \beta \circ F \circ \alpha $ in $\S$ has a \bpt at $\one$ 
and
\newline
$\phi(\one) \neq 1 $.

\item[\rm(iv)]
 $\lim_{y \to \i} F(iy\one) = C \inn \R$ and
$$
 \liminf_{y \to \i} \ y \| v_{(iy\one)} \|\, < \, \i .
$$

\item[\rm(v)]
The vector $v$ is in the domain of $X$.
\end{enumerate}
\et

%% file: matmonSep13.tex
\section{Locally matrix monotone functions}
\label{secm}
Recall the definition of locally $n$-matrix monotone.

{\bf Definition~\ref{defi10}\ }{\em
Let $E$ be an open set in $\R^d$, and $f$ be a real-valued
$C^1$ function on
$E$.
Say $f$ is locally $\mn$-monotone on $E$ if, whenever
$S$ is 
in $\csmn$
with $\sigma(S) = \{ x_1,\dots, x_n \}$
consisting of $n$ distinct points in  $E$, and
$S(t)$ is a $C^1$ curve 
in $\csmn$
with $S(0) = S$ and
$\dis \frac{d}{dt} S(t) |_{t=0} \geq 0 $,
then
$\dis \frac{d}{dt} f(S(t)) |_{t=0}$
exists and is $ \geq 0 $.
}

If $S$ is in $\csmn$,
we can choose an orthonormal basis of eigenvectors  
that diagonalize all the $S^r$'s simultaneously, so 
\be
\label{eqm01}
S^r \= \left(
\begin{array}{ccc}
x_1^r && \\
&\ddots& \\
&& x_n^r
\end{array}
\right) \qquad \forall \ 1 \leq r \leq d .
\ee
If $S(t)$ is a $C^1$ curve of commuting self-adjoints, then 
$ S(0) + t \,  \frac{d}{dt} f(S(t)) |_{t=0}$ commutes to first order.

For any $X \in M_n$ we define $\diag X $ to be the diagonal matrix in $M_n$ with diagonal entries $X_{ii}$, and for any $\Delta\in SAM_n^d$ we define $\diag \Delta$ to be $(\diag\Delta^1, \dots, \diag\Delta^d)$.

\bd
We shall say that $S$ in $\csmn$ is generic if
its spectrum consists of $n$ distinct points.
\ed

\bl
\label{lemma}
Let $S$ be in $\csmn$ and $\Delta$ be in $\smn$, with  $S$ generic.
Then there exists a $C^1$ curve $S(t)$ of commuting self-adjoints
with $S(0) = S$ and
$S'(0) = \Delta$ if and only if 
\be
\label{eqm1}
[ S^r, \Delta^s] \= [ S^s , \Delta^r] \qquad \forall \ 1 \leq r \neq s \leq d .
\ee
\el
\bp
($\Rightarrow$): If $S(t) = S + t \De + o(t)$ is commutative, calculate
$$
[S^r (t) , S^s(t)] \= t \left( [ S^r, \Delta^s] - [ S^s , \Delta^r] \right) + o(t).
$$
The coefficient of $t$  must vanish, giving (\ref{eqm1}).

($\Leftarrow$):
Suppose  $S$ is as in (\ref{eqm01}),
and (\ref{eqm1}) holds. This means
\be
\label{eqm04}
\De^s_{ij} ( x^r_j - x^r_i ) \=
\De^r_{ij} ( x^s_j - x^s_i ) \qquad \forall \, r \neq s ,
\ee
so
\be
\label{eqm06}
\De^r_{ij} \frac{1}{x^r_j - x^r_i} \=
\De^s_{ij} \frac{1}{x^s_j - x^s_i}  \qquad {\rm if\ }
x^r_j - x^r_i \neq 0 \neq
x^s_j - x^s_i.
\ee
Define a skew-symmetric matrix $Y$ by
\be
\label{eqm07}
Y_{ij} \= 
\De^r_{ij} \frac{1}{x^r_j - x^r_i} \qquad {\rm for\ any\ }r\ {\rm such\ that\ }
x^r_j - x^r_i \neq 0 .
\ee
For any $i \neq j$, there is some $r$ with 
$
x^r_j - x^r_i \neq 0$, so (\ref{eqm07}) defines $Y_{ij}$;
and (\ref{eqm06})
says it doesn't matter which $r$ we choose.
Let all the diagonal terms of $Y$ be $0$.

Define
\be
\label{eqz3}
S^r(t) \= e^{tY}( S^r +  t\ \diag \Delta^r)  e^{-tY} .
\ee
Since $e^{tY}$ is a unitary matrix and $S^r  + t\ \diag \Delta^r$ is diagonal, $S(t) \in CSAM_n^d$ and
\[
 \frac{d}{dt}S^r(t)|_{t=0} = [Y,S^r] + \diag \Delta^r = \Delta^r.
\]
\ep


If $S$ and $\Delta$ satisfy (\ref{eqm1}) and $S$ is generic
then for any function $f$ that is $C^1$ on a neighborhood of $\sigma(S)$ 
we define the directional derivative of $f$ at $S$ in direction $\Delta$ by
\be
\label{eqx3}
D_\Delta f(S) = \frac{d}{dt}f(S(t))|_{t=0}
\ee
where $S(t)$ is the curve  given by equations (\ref{eqz3}) and (\ref{eqm07}).  
We shall show in Proposition~\ref{prx4} that 
(\ref{eqx3}) is actually unchanged if $S(t)$ is replaced by any other curve
that agrees with it to first order. First, let us show that the right-hand side
of (\ref{eqx3}) exists.
Indeed, 
\be\label{formfSt}
 f(S(t)) = e^{tY} f(S+t\ \diag\Delta) e^{-tY}.
\ee
 Since $S+t\ \diag\Delta$ is diagonal, $f(S+t\ \diag\Delta)$ is diagonal, with $i$th entry 
\[
 f(x_i + t \Delta_{ii}) = f(x_i) + t\sum_{r=1}^d \Delta_{ii}^r \frac{\partial f}{\partial x^r} (x_i) + o(t).
\]
In other words,
\[
f(S+t\ \diag\Delta) = f(S) + t \sum_{r=1}^d (\diag \Delta^r)\frac{\partial f}{\partial x^r}(S) +o(t).
\]
Hence, on differentiating equation (\ref{formfSt}) at $0$ we obtain
\[
\frac{d}{dt}f(S(t))|_{t=0} = [Y,f(S)] + \sum_{r=1}^d (\diag \Delta^r)\frac{\partial f}{\partial x^r}(S).
\]
We have shown the following.
\begin{prop}\label{7.8}
Let $S$ 
be a generic $d$-tuple of commuting self-adjoint matrices in $M_n$. Fix an orthonormal
basis 
of eigenvectors, so every $S^r$ is diagonal: 
$$
S^r \= \left( 
\begin{array}{ccc}
x_1^r && \\
&\ddots& \\
&& x^r_n 
\end{array} \right) .
$$
Let $\Delta$ be a $d$-tuple of self-adjoints satisfying (\ref{eqm1}).
Let $f$ be $C^1$ on a neighborhood of $\{x_1,\dots,x_n \}$ in $\R^d$,
where each $x_j$ is the $d$-tuple $(x^1_j, \dots, x^d_j)$.
Then
\be
\label{eqma2}
\left[ D_\De f (S) \right]_{ij}
\= \left\{
\begin{array}{ll}
\De^r_{ij} \frac{f(x_j) - f(x_i)}{x_j^r - x_i^r} &
\mbox{if\ }i\neq j,\ \mbox{where\ }x_j^r \neq x_i^r  \\
&\\
\sum_{r=1}^d 
\De^r_{ii} \frac{\partial f}{\partial x^r} |_{x_i} & \mbox{if \ } i = j.
\end{array}
\right.
\ee   
\end{prop}
%

\begin{cor} \label{DDelfg}
For $S, \Delta$ as in Proposition \ref{7.8}, if $f, g$ are $C^1$ functions that agree to first order on $\sigma(S)$, then $D_\Delta f(S) = D_\Delta g(S)$.
\end{cor}

\begin{lemma}\label{sep}
Let $R$ and $S$ be in $CSAM_n^d$. For every point $\mu$ in the joint spectrum of $R$
there is an $x_p$ in the joint spectrum of $S$ with 
\be
\label{eqsep}
 \| \mu - x_p \|\  \leq\   \sqrt{dn} 
\| R - S \| .
\ee
\end{lemma}
\begin{proof}
Choose an orthonormal basis that diagonalizes $S$, so that $S$ is as in (\ref{eqm01}).
Let $\De = R - S$.
Let $\mu$ be a joint eigenvalue of $R$ with corresponding eigenvector
$\xi = (\xi_1, \dots, \xi_n)^t$. Choose $p$ so that $|\xi_p | \geq |\xi_j |$ for all $1 \leq j \leq n$.

Then for each $1 \leq r \leq d$, we have $$
R^r \xi \= \mu^r \xi \= (S^r + \De^r) \xi .
$$
So in particular,
$$
\sum_{j=1}^n R^r_{pj} \xi_j \= \mu^r \xi_p .
$$
Therefore
\[
(\mu^r - x_p^r  ) \xi_p \= \sum_{j =1}^n \De^r_{pj} \xi_j .
\]
So
\beq
| \mu^r - x_p^r | & \ \leq \ & \sum_{j=1}^n | \De^r_{pj} | \\
&\leq & \sqrt{n} \sqrt{ \sum_j |\De^r_{pj} |^2 } \\
&\leq & \sqrt{n} \| \De^r \| ,
\eeq
and hence 
$$
\sum_{r=1}^d
| \mu^r - x_p^r |^2 \ \leq \
dn \| \De \|^2.
$$
%
\end{proof}

\begin{lemma}\label{rellich}
If  $R(t)$ is a Lipschitz path in $CSAM_n^d, \ 0\leq t <1,$ with $R(0) = S$ generic
then
 there exists $\vare > 0 $
and Lipschitz maps
 $X_1, \dots , X_n:[0,\vare)\to \R^d$  such that $\sigma(R(t)) =\{X_j(t): j=1,\dots,n\}$.  
\end{lemma}
\begin{proof}  

Choose an orthonormal basis that diagonalizes $S$, so $S$ is as in (\ref{eqm01}).
The joint eigenvalues of $S$ are the points $x_i = (x_i^1,\dots,x_i^d)$, and
genericity means $\| x_i - x_j \| > 0$ if $i \neq j$.
Choose $\vare$ so that for all $0 \leq t \leq \vare$,
\be \label{eqz5}
\sqrt{dn} \| R(t) - S  \| \ \leq \ \frac{1}{3} \min_{i\neq j}\  \|x_i - x_j \|.
\ee
%

By Lemma~\ref{sep}, for every joint eigenvalue $x$ of $S$ there is a joint eigenvalue $\mu$ of
$R(t)$ within $\sqrt{dn} \| R(t) - S \|$ of it. By (\ref{eqz5}), this means that $R(t)$ is also generic,
and  it makes sense to talk  of the 
joint eigenvalue of $R(t)$ that is closest to  $x_j$.
Let us call these joint eigenvalues $X_j(t)$. We have proved that 
\[
\| X_j (t) - x_j \| \ \leq \ \sqrt{dn} \| R(t) - S \|
\qquad \forall \ 0 \leq t \leq \vare .
\]
Repeating the argument with $R(t_1)$ in place of $S$, we get 
\[
\| X_j (t_2) - X_j (t_1) \| \ \leq \ \sqrt{dn} \| R(t_2) - R(t_1) \|
\qquad \forall \ 0 \leq t_1, t_2 \leq \vare .
\]
As $R$ is assumed to be Lipschitz, we get that each $X_j$ is Lipschitz also.
\end{proof}

\begin{prop}
\label{prx4}
If $S$ is generic in $CSAM_n^d$, $\Delta$ is in $SAM^n_d$, and they 
 satisfy the commutation relations (\ref{eqm1}),
 then for any $C^1$ path $R(t) \in CSAM_n^d$ such that $R(0)=S, \ R'(0)=\Delta$ and any $f \in C^1$,
\be
\label{propDDelta}
 \frac{d}{dt}f(R(t))|_{t=0} = D_\Delta f(S).
\ee
\end{prop}
\begin{proof}
If $g$ is a monomial then a simple calculation shows that
\[
\frac{d}{dt}g(R(t))|_{t=0}
\]
exists and depends only on $g, S$ and $\Delta$.  It follows that, for any polynomial $g$,
\be \label{3.1}
\frac{d}{dt}g(R(t))|_{t=0} = \frac{d}{dt}g(S(t))|_{t=0} = D_\Delta g(S).
\ee

Consider any $f\in C^1$ and pick a polynomial $g$ that agrees with $f$ to first order on $\sigma(S)$.  By Corollary \ref{DDelfg},
\be\label{fgS}
D_\Delta f(S) = D_\Delta g(S).
\ee
We claim that
\be\label{dfgR}
\frac{d}{dt}g(R(t))|_{t=0} = \frac{d}{dt} f(R(t))|_{t=0}.
\ee
For by Lemma \ref{rellich} there exist Lipschitz functions $X_1, \dots, X_n: [0,\vare) \to \R^d$ such that 
$\sigma(R(t)) = \{X_1(t), \dots, X_n(t)\}$ for all $t$.  Then $f(S)=g(S)$ and
\beq
\|(f-g)(R(t))\| &\= & \max_i |(f-g)(X_i(t))| \\
&=&  o( \max_i \|X_i(t) - X_i(0) \| ) \\
&=& o(t).
\eeq
Hence
\[
\left\| \frac{f(R(t))-f(S)}{t} -\frac{g(R(t)) - g(S)}{t} \right\| \to 0 \mbox{ as } t\to 0.
\]
In view of equation (\ref{3.1}),
\[
\frac{f(R(t))-f(S)}{t} \to \frac{d}{dt}g(R(t))|_{t=0} =D_\Delta g(S) \mbox{ as } t\to 0.
\]
On combining this relation with equation (\ref{fgS}) we obtain equation (\ref{propDDelta}).
\end{proof}

\begin{cor}
A real-valued $C^1$ function $f$ on an open set $E \subseteq \R^d$
is locally $M_n$-monotone if and only if 
\[
D_\Delta f(S) \ \geq \ 0 
\]
for every generic $S$ in $CSAM_n^d$ with spectrum in $E$ and every
$\Delta$ in $SAM^d_n$ such that $\Delta \geq 0 $ and 
\[
[ S^r, \Delta^s] \= [ S^s , \Delta^r] \qquad \forall \ 1 \leq r \neq s \leq d .
\]
\end{cor}
The statement follows immediately from Definition~\ref{defi10} and Proposition \ref{prx4}.

We can now characterize locally matrix monotone functions.

\bt
\label{thmma}
Let $E$ be an open set in $\R^d$, and $f$ a real-valued $C^1$ function on $E$.
Then $f$ is locally $\mn$-monotone if and only if $f$ is in $\lnde$.
\et
\bp
($\Leftarrow$)
We must show: if $S$ is generic with $\sigma(S) \subset E$, if $\De$ is 
a positive $d$-tuple and $[S^r,\De^s] = [S^s,\De^r]$ for all $r,s$, then
$ D_\De f(S) \geq 0 $.

Let $\sigma(S) = \{ x_1,\dots,x_n \}$.
Choose $A^r$ as in  Definition~\ref{deb1}. 
For $i \neq j$, assume without loss of generality that $x^1_j \neq x^1_i$. Then 
\beq
[ D_\De f(S) ]_{ij} &\=& 
\De^1_{ij} \frac{f(x_j) - f(x_i)}{x_j^1 - x_i^1} \\
&=& 
\frac{\De^1_{ij}}{x^1_j - x^1_i} \left( \sum_{r=1}^d (x_j^r - x_i^r) A^r(i,j) \right) \\
&=&
\sum_{r=1}^d \De^r_{ij} \, A^r(i,j) .
\eeq
(We get the last line by using (\ref{eqm06})).
By (\ref{eqma2}) the same formula holds 
for $ 
[ D_\De f(S) ]_{ij} $
when $i=j$, so
$
D_\De f(S) $
is the sum of the Schur products of $\De^r$ with $A^r$, so is positive.

($\Rightarrow$)
Let 
$f$ be locally $\mn$-monotone, and fix $\{ x_1, \dots, x_n \}$ distinct
points in $E$. Let $S$ be given by (\ref{eqm01}).
We wish to find positive matrices $A^r$
such that (\ref{defb12}) and (\ref{defb11}) hold.

Let $\mathcal G$ be the set of all skew-symmetric real
$n$-by-$n$ matrices
$\Gamma$
with the property that there exists a $d$-tuple $A$ of real
positive semi-definite
matrices satisfying
\se\att
\begin{eqnarray}
\label{eqm0a1}
A^r(i,i) &=& \left. \frac{\partial f}{\partial x^r} \right|_{x_i}
\quad 1\leq i \leq n,\ 1 \leq r \leq d
\\
\sum_{r=1}^d (x^r_j - x^r_i) \, A^r(i,j)
&=&
\Gamma_{ij}
\qquad 1 \leq i \neq j \leq n.
\att
\label{eqm0a2}
\end{eqnarray}
Let $\La$ be the matrix $\La_{ij} = f(x_j) - f(x_i)$. We wish to show $\La$ is in $\mathcal G$.

Notice that $\mathcal G$ is a closed convex set. Moreover, it is non-empty, because
$\dis \left. \frac{\partial f}{\partial x^r} \right|_{x_i}  $ 
is always greater than
or equal to $0$, so $\mathcal G$ contains $0$.
Without loss of generality, we shall assume that 
\[
f_{r,i} \ := \ \left. \frac{\partial f}{\partial x^r} \right|_{x_i} \ > \ 0, \quad
\forall\ 1 \leq r \leq d, \ 1 \leq i \leq n .
\]
(This can be done by adding $\epsilon ( x^1 + \dots + x^d)$ to $f$,
and then letting $\epsilon \to 0$ at the end of the argument.)

If $\La$ is not in $\mathcal G$, by the Hahn-Banach theorem
there is a real  linear functional $L :\mn \to \R$ and some $\delta \geq 0$
so that $L(\Gamma) \geq - \delta$ for all $\Gamma $ in
$\mathcal G$,  and $L(\La) < -\delta$.
Any such linear functional is of the form $L (T) = {\rm tr} (TK)$ for some real  matrix $K$. 
Replacing $K$ by $\frac{1}{2} (K -  K^t)$ will not change the value of $ L$ on
skew-symmetric real matrices, so we can assume that there is a real skew-symmetric
matrix $K$ such that
${\rm tr} (\Gamma K) \geq  - \delta$ for all $\Gamma $ in $\mathcal G$,
and 
${\rm tr} (\La K) < - \delta$.

Define $\De$ by 
\[
\label{eqmnew2}
\De^r_{ij} \= 
(x^r_j - x^r_i ) K_{ji} , \quad i \neq j,
\]
and with the diagonal entries $\Delta^r_{ii}$  chosen so that
each $\De^r \geq 0$   and so that 
\be
\label{eqmnew4}
\mu^r \ := \
\sum_{i=1}^n f_{r,i} \De^r_{ii}
\ee
is minimal over all choices 
of $\Delta^r_{11}, \dots, \Delta^r_{nn}$
such that $\De \geq 0$. (A minimal choice exists, since all the $f_{r,i}$ are strictly positive by assumption).
Then $\De$ is in $\smn$, and 
$$
[ \De^s, S^r ]_{ij} \= (x_j^s - x_i^s) K_{ji} (x_j^r - x^r_i) \= [\De^r, S^s ]_{ij},
$$
so $\De$ satisfies (\ref{eqm1}). 

As $f$ is locally $\mn$ monotone, 
we must have then that $D_\De f(S) \geq 0$ by Lemma 7.3.

By assumption, 
\begin{eqnarray}
\nonumber
- \delta \ > \ {\rm tr}(\La K) &\=&
\sum_{1 \leq i \neq j \leq n} [f(x_j) - f(x_i) ] K_{ji} \\
\label{eqxj3}
&=&
\sum_{1 \leq i \neq j \leq n} \frac{f(x_j) - f(x_i)}{x^r_j - x^r_i} \De^r_{ij} \\
\nonumber
&=& \sum_{1 \leq  i , j \leq n}
[ D_\De f(S) ]_{ij}  - \sum_{r=1}^d \sum_{i=1}^n \De^r_{ii} f_{r,i},
\end{eqnarray}
\att
where in \eqref{eqxj3}, for each $(i,j)$ we choose a (perhaps different) $r$ so that the
denominator is non-zero.
Therefore
\be
\label{eqmnewy}
\sum_{r=1}^d \mu^r - \delta \ > \  \sum_{1 \leq  i , j \leq n}
[ D_\De f(S) ]_{ij} \geq 0.
\ee
So if we can prove that
\be
\label{eqmnew1}
\sum_{r=1}^d \mu^r \ \leq \ \delta,
\ee
we shall derive a contradiction.

By Duffin's strong duality theorem \cite{duf56}, the minimun $\mu^r$ in \eqref{eqmnew4}
satisfies
\be\label{mur}
-\mu^r \ = \ \min \sum_{i\neq j} \Delta_{ij} A^r(i,j),
\ee
where $A^r$ range over the set of
real positive matrices  such that the diagonal entries of $A^r$ are $f_{r1}, \dots,f_{rn}$ for each $r$.

For each such $A = (A^1, \dots, A^d)$,  let $\Gamma$ be the corresponding element of $\mathcal G$:  $\Gamma_{ii}=0$ and
\[
\Gamma_{ij}= \sum_{r=1}^d (x^r_j-x^r_i) A^r(i,j) \quad \mbox{ for }i\neq j.
\]
We have
\begin{align*}
-\delta &\leq \mathrm{tr}~\Gamma K \\
	&= \sum_{i\neq j} \sum_{r=1}^d  (x^r_j - x^r_i) A^r(i,j) K_{ji} \\
	&= \sum_{r=1}^d  \sum_{i\neq j} \Delta^r_{ij} A^r(i,j).
\end{align*}
Hence, by equation \eqref{mur}, $-\delta \leq \sum_{r=1}^d (-\mu^r)$,
so $\sum_{r=1}^d \mu^r \ \leq \ \delta$. This contradicts
\eqref{eqmnewy}, so
it follows that $\La \in\mathcal G$, and necessity is proved.

\ep

As the dimension of  the matrices increases, the condition that a 
function $f$ be locally monotone becomes more stringent. On an infinite dimensional
Hilbert space, the requirement becomes that $f$ 
be in the Loewner class, as we shall see in the next section.

%
%
%
%

%% file: locopmon.tex
\section{Locally operator monotone functions}
\label{seclo}
We defined locally operator monotone functions in Definition~\ref{defilocop}.
We shall show that being locally operator monotone is the same as 
being locally $M_n$-monotone for all $n$, which in turn is the same as
being in the L\"owner class $\L(E)$.
\bt
\label{thmloa}
Let $E$ be an open set in $\R^d$, and $f$ a real-valued $C^1$ function on $E$.
The following are equivalent.

(i) The function $f$ is locally $\mn$-monotone on $E$ for all $ n \geq 1$.

(ii) The function $f$ is in $\L(E)$.

(iii) The function $f$ is locally operator monotone on $E$.
\et
\vs
The equivalence of $(i)$ and $(ii)$ follows from Theorem~\ref{thmma} and
Proposition~\ref{propba1}. The implication $ (iii) \Rightarrow (i)$ is obvious.
We need to prove that $(ii) \Rightarrow (iii)$.
First we need some preliminary results.

\bprop
\label{proploa}
Let $E$ be an open set in $\R^d$, and let $f \, \in \, \L(E)$.
Then there is a model $(\M,v)$ for $f$
such that $v_z$ is locally Lipschitz on $E$.
\eprop
\bp
By Proposition~\ref{propba1} 
we can 
extend $f$ to a function $F$ in $\L$
that extends analytically 
across $E$ and agrees with $f$ on $E$.
For this $F$ we have a model $(\M,v)$
so that
\be
\label{eqlo2}
F(z) -   \overline{F(w)} \= \la (z-w^*) v_z, v_w \ra_{\M}
\qquad \forall \ z,w \, \in \, E \cup \Pid ,
\ee
and by Proposition~\ref{propmodc},
if $w$ is in $E$, then
$v_w$ is the limit of $v_z$ as $z$ tends to $w$ non-tangentially from inside $\Pid$.

Fix $w$ in $E$ (so $F(w)$ is real). Then, by analyticity, we have
for $z$ close to $w$ :
\be
\label{eqlo3}
F(z) - F(w) \= \sum_{r=1}^d  
\left. \frac{\partial f}{\partial x^r} \right|_{w}
(z^r - w^r) \ + \
\sum_{|\a| \geq 2} 
\left. \frac{\partial^\a f}{\partial x^\a} \right|_{w}
\frac{(z-w)^\a}{\a!} 
. \ee
From (\ref{eqlo2}), we get
\beq
F(z) -   {F(w)} &\=& \la (z-w) v_z, v_w \ra
\\
&\=&
\sum_{r=1}^d (z^r - w^r) \la v^r_w, v^r_w \ra
\ + \
\sum_{r=1}^d (z^r - w^r) \la v^r_z - v^r_w, v^r_w \ra
. \eeq
As $z$ tends to $w$ non-tangentially, the second term
is $o(\| z - w \|)$, so comparing with (\ref{eqlo3}) we see that
\be
\label{eqlo4}
\| v_w^r \|^2 \= 
\left. \frac{\partial f}{\partial x^r} \right|_{w}
\qquad \forall \, 1 \leq r \leq d .
\ee
Now let $z$ and $w$ both be in $E$. Comparing (\ref{eqlo2}) and
(\ref{eqlo3}), we get
\be
\label{eqlo5}
\la (z-w) v_z, v_w \ra
\ - \ \la (z-w) v_w, v_w \ra
\=
\sum_{|\a| \geq 2} 
\left. \frac{\partial^\a f}{\partial x^\a} \right|_{w}
\frac{(z-w)^\a}{\a!} 
. \ee
Swapping $z$ and $w$, we get
\be
\label{eqlo6}
\la (w-z) (v_w - v_z), v_z \ra
\=
\sum_{|\a| \geq 2} 
\left. \frac{\partial^\a f}{\partial x^\a} \right|_{z}
\frac{(w-z)^\a}{\a!} 
. \ee
Subtracting (\ref{eqlo5}) from (\ref{eqlo6}), we get
\be
\label{eqlo7}
\la (z-w) (v_z - v_w), (v_z - v_w) \ra
\=
\sum_{|\a| \geq 2} 
\left(
(-1)^{|\alpha|} \left. \frac{\partial^\a f}{\partial x^\a} \right|_{z}
\ -\ 
\left. \frac{\partial^\a f}{\partial x^\a} \right|_{w}
\right)
\frac{(z-w)^\a}{\a!} 
. \ee
But since $f$ is analytic,
\[
\left(
\left. \frac{\partial^\a f}{\partial x^\a} \right|_{z}
\ -\ 
\left. \frac{\partial^\a f}{\partial x^\a} \right|_{w}
\right)
\= O(\| z - w \|) ,
\]
and so the right-hand side of (\ref{eqlo7}) is $O(\| z -w \|^3)$.
Therefore
\be
\label{eqlo8}
\la (z-w) (v_z - v_w), (v_z - v_w) \ra
\=
O(\| z -w \|^3) .
\ee
If all the differences 
$|z^r - w^r|$ were comparable, we could conclude immediately that
$\| v_z - v_w \| = O(\|z-w\|)$.
If they are not, we can get round this difficulty by connecting $z$ to $w$
by two line segments.

Indeed, suppose $\max_{1 \leq r \leq d} |z^r - w^r | = \vare$.
Choose numbers $a^r$ and $b^r$ with modulus between $1/2$ and $2$
so that 
$$
z^r - w^r \= (a^r - b^r) \vare \qquad \forall \ 1 \leq r \leq d .
$$
Let
$$
x^r \= w^r + a^r \vare \= z^r + b^r \vare.
$$
Then applying (\ref{eqlo8})
to the pairs $(z,x)$
and $(x,w)$, we get 
\beq
\| v_z - v_w \| &\ \leq \ &
\| v_z - v_x \| + 
\| v_x - v_w \| \\
&=&
O(\| z - x \| + \| x - w \|) \\
&=&
O(\| z - w \|) ,
\eeq
as desired.
\ep

%
%
%
%

\bs
Suppose now $E, f $ and $(\M,v)$
are as in Proposition~\ref{proploa}.
So $v : z \mapsto v_z$ is a map from $E$ to $\M$. Let $S $ be a $d$-tuple
of bounded commuting self-adjoint operators on a Hilbert space $\h$, with $\sigma(S)
\subset E$. We want to define an operator $\tilde v(S) \, \in \, 
B(\h, \h \otimes \M)$.

We do this by choosing an orthonormal basis for $\M$, and writing
\be
\label{eqar3}
v(z) \ :=\ v_z \= 
\left(
\begin{array}{c}
v_1(z) \\
v_2(z) \\
\vdots
\end{array}
\right) .
\ee
We must caution the reader that the subscripts on the right-hand side of (\ref{eqar3})
run over the dimension of $\M$, and are not interchangeable with the superscripts on $v$
that identify which piece of $\M^1, \dots, \M^d$ one is in.
After using the orthonormal basis to identify $\M$ with $\ell^2$, we have
\[
v_z \=
\left[
\begin{array}{c}
v_z^1 \\
\vdots \\
v_z^d
\end{array}
\right]
\=
\left(
\begin{array}{c}
v_1(z) \\
v_2(z) \\
\vdots
\end{array}
\right) .
\]
Then
\be
\label{eqlo10}
\tilde{v}(S) \ := \
\left(
\begin{array}{c}
v_1(S) \\
v_2(S) \\
\vdots
\end{array}
\right) \ :\, \h \to \h \otimes \M.
\ee
The operator $\tilde{v}(S)$ is bounded, because if $S$ has spectral
measure $\Lambda$ and $\xi$ is a unit vector in $\h$, then
\se\att
\begin{eqnarray}
\label{eqlo56}
\| \tilde{v}(S) \xi \|^2 &\=&
\sum_j \int_{\sigma(S)} |v_j |^2 d\langle \Lambda \xi,\xi \ra
\\
\nonumber
&=&
\int_{\sigma(S)} \sum_{r=1}^d \frac{\partial f}{\partial x^r} d \langle \Lambda \xi,\xi \ra
\\
\att
\label{eqx7}
&\leq&
\sup_{z \, \in \, \sigma(S)} \
\left. \sum_{r=1}^d \frac{\partial f}{\partial x^r} \right|_z ,
\end{eqnarray}
and the last sum is finite because $\sigma(S)$ is compact
and $f$ is $C^1$.

The operator $\tilde v(S)$ does not depend on the choice of orthonormal
basis in $\M$.  A simple calculation shows that, for any $h\in\mathcal{H}$
and $m\in\M$,
\[
\tilde v(S)^* (h\otimes m) = \langle m, v( \cdot )\rangle (S) h.
\]
This gives a coordinate-free expression for $\tilde v(S)^*$, hence also
for $\tilde v(S)$.


\bl
\label{lemlo1}
 Let $E, \ f$ and $(\M, v)$ be as in Proposition \ref{proploa},
and $S $ be a $d$-tuple
of bounded commuting self-adjoint operators on a Hilbert space $\h$, with $\sigma(S)
\subset E$.  Then
\be
\label{eqx8}
  \|\tilde v(S)\| \  \leq  \ \left(\sup_{\sigma(S)} \sum_{r=1}^d \frac{\partial f}{\partial x^r}
\right)^{\tfrac 12} .
\ee
Moreover, $\tilde v$ is continuous.
\el
\bp
Inequality (\ref{eqx8}) has been proved in (\ref{eqx7}).
To prove continuity of $\tilde v$, 
let $K$ be a compact
subset of $E$ with $\sigma (S) \subset int(K) \subset E$.
Let $\vare > 0$.

As $v$ is continuous on $K$ and $K$ is compact,
there exists $N$ such that
$\sum_{j =N+1}^\i |v_j(z)|^2 \leq \vare^2/9 $ on $K$.
For $1 \leq j \leq N$, there is a polynomial $p_j$ such that
$\| p_j - v_j \|_{\i} \leq \vare/9N$ on $K$.
There exists $\delta > 0$ so that if $\|T^r - S^r \| \leq \delta$, then
$\sigma(T) \subseteq K$ and
$\| p_j (T) - p_j(S) \| \leq \vare/9N$ for each $1 \leq j \leq N$.

Let $\tilde{v}_N (S)$ be the operator
$$
\left(
\begin{array}{c}
v_1(S) \\
\vdots \\
v_N(S) \\
0 \\
\vdots
\end{array} \right) .
$$
Then $\| \tilde{v}_N(S) - \tilde{v}(S) \| \leq \vare/3$ by (\ref{eqlo56}), and
similarly
$\| \tilde{v}_N(T) - \tilde{v}(T) \| \leq \vare/3$.
As
$$
\| v_j (T) - v_j(S) \| \ \leq \ 
\| v_j (T) - p_j(T) \| + \| p_j(T) - p_j(S) \| +
\| p_j(S) - v_j (S) \| ,
$$
and each entry is at most $\vare/9N$, 
$$
 \|\tilde v_N(S) - \tilde v_N(T) \| \leq N\left( \frac{\vare}{9N} + \frac{\vare}{9N} + \frac{\vare}{9N}
  \right) = \frac{\vare}{3},$$
and hence
$$
\| \tilde{v}(T) - \tilde{v}(S) \| \ \leq \ \vare .
$$
\ep

\bs
We assume that $\M$ is decomposed as 
$\M = \M^1 \oplus \dots \oplus \M^d$, and
$P^r$ is the orthogonal projection from $\M$ onto $\M^r$.
If $S = (S^1, \dots, S^d)$ is a $d$-tuple of operators on $\h$, 
we shall write
\be
\label{eqlo35}
S \odot I \ := \ S^1 \otimes P^1 \oplus \dots \oplus S^d \otimes P^d
\ee
which is an operator on $\h \otimes \M$.

\bprop
\label{proplob}
Let $E$ be open in $\R^d$, let $f \, \in \, \L(E)$, and assume $(\M,v)$ is a 
model of $f$ 
for which $v$ is continuous.
Let $S$ and $T$ be $d$-tuples of commuting self-adjoint operators 
on a Hilbert space $\h$ with spectrum in 
$E$.
Then
\be
\label{eqlo11}
f(T) - f(S) \=
\tilde{v}(S)^* \left[ T \odot I - S \odot I \right] \tilde{v}(T) .
\ee
\eprop
\bp
First assume that $S$ and $T$ are (separately) diagonalizable.
Let $\xi$ be an eigenvector of $S$, and $\eta$ an eigenvector of $T$, 
so for some numbers $z^r,w^r$ we have
\beq
S^r \xi &\=& w^r \xi \\
T^r \eta &=& z^r \eta \qquad \forall \ 1 \leq r \leq d.
\eeq
Then
\[
\la \left[ f(T) - f(S) \right] \eta , \xi \ra_\h \=
\la \left[ f(z) - f(w)^* \right] \eta , \xi \ra_\h .
\]
Also,
\beq
\lefteqn{
\la 
\tilde{v}(S)^* \left[ T \odot I - S \odot I \right] \tilde{v}(T) 
 \eta , \xi \ra_\h 
} \\
&{} \quad = \ &
\la
\left[ T \odot I - S \odot I \right] \eta \otimes v(z) ,
\xi \otimes v(w) \ra_{\h \otimes \M}
\\
&{} \quad = \ &
\sum_{r=1}^d
\la \eta, \xi \ra_\h \ 
\la (z^r - \bar{w^r}) v^r(z) , v^r(w) \ra_{\M^r}
\\
&{} \quad = \ &
(f(z) - f(w) ) 
\la \eta, \xi \ra_\h .
\eeq
So both sides of (\ref{eqlo11}) agree if you apply them to an
eigenvector of $T$ and then take the inner product with an eigenvector
of $S$.
By linearity, this is true also for linear combinations
of eigenvectors, and as these are assumed dense, we get that
(\ref{eqlo11}) holds.

If $S$ and $T$ are not diagonalizable, by the spectral theorem
we can approximate them in norm by operators that are, and as $\tilde v$ and
$f$ are both continuous, one gets (\ref{eqlo11}) in the limit.
\ep

{\sc Proof of Theorem~\ref{thmloa}:}
Assume $f$ is in $\L(E)$, and $S(t)$ is a curve of commuting self-adjoint
$d$-tuples with $S(0) = S$ and $S'(0) = \Delta \geq 0$.
Choose a model $(\M,v)$ with $v$ locally Lipschitz. Then by Proposition~\ref{proplob},
\[
f(S(t)) - f(S) \= 
\tilde{v}(S)^*
\left[ (S(t) -S) \odot I \right] \tilde{v}(S(t)) .
\]
As
\[
S(t) = S + t \Delta + o(t) , \]
we get
\beq
\left. \frac{d}{dt} f(S(t)) \right|_{0}
&\=  &
\lim_{t \to 0}
\tilde{v}(S)^*
\left[ \Delta \odot I \right] \tilde{v}(S(t)) 
\ + \
\lim_{t \to 0}
\tilde{v}(S)^*
\left[ o(1) \right] \tilde{v}(S(t)) 
\\
&=&
\tilde{v}(S)^*
\left[ \Delta \odot I \right] \tilde{v}(S) .
\eeq
Hence $f(S(t))$ is differentiable at $0$, and its derivative is a positive
semi-definite operator.
\ep

%% file: glob.tex
\section{Globally Operator Monotone Functions}
\label{secg}

\bd
Let $E$ be an open set in $\R^d$, and $f$ be a real-valued
$C^1$ function on
$E$.
Say $f$ is globally operator monotone on $E$ if, whenever
$S$ and $T$
are $d$-tuples of commuting bounded self-adjoint operators on a Hilbert space
with $\sigma(S) \cup \sigma(T) \subset E$,
and
$S \leq T$, then
$f(S) \leq f(T)$.
\ed

If $F$ has the form in Theorem~\ref{thmd3}, then it
is globally monotone on boxes in the $\mu$-resolvent of $X$.

\bt
\label{thmg1}
Let $X$ be a densely-defined self-adjoint operator on a graded Hilbert space $\M
\= \M^1 \oplus \dots \oplus \M^d$, let $v \in \M$, and let $F$ be given by
\[
F(z) \= \la (X-z)^{-1}v,v \ra .
\]
Let $E$ be an open box in $\R^d$ that is in the $\mu$-resolvent of $X$.
Then $F$ is globally operator monotone on $E$.
\et
\bp
First observe that 
if $S$ is a commuting $d$-tuple of self-adjoint operators on $\h$ and
$\sigma(S) \subset E$, then
\be\label{formFS}
F(S) = R_v^* (I_\mathcal{H} \otimes X - S \odot I)^{-1} R_v
\ee
where $\odot$ is as in equation (\ref{eqlo35}) and
\begin{eqnarray*}
 R_v \ :\  \mathcal { H} &\to & \mathcal{H \otimes M} \\
 h &\mapsto & h\otimes v.
\end{eqnarray*}
Thus equation (\ref{formFS}) means that 
for any vectors $\xi $ and $\eta$ in $\h$,
\be
\label{eqg2}
\la F(S) \eta , \xi \ra_\h
\=
\la
( I_{\h} \otimes X - \sum_{r=1}^d S^r \otimes P^r )^{-1} \ \eta \otimes v, 
\xi \otimes v \ra_{\h \otimes \M} . 
\ee
Indeed, if $\eta$ is an eigenvector of $S$ with eigenvalues $a^r$,
then $F(S) \eta = F(a) \eta$, 
so the left-hand side of (\ref{eqg2}) is $F(a) \la \eta , \xi \ra$.
But we have\[
( I_{\h} \otimes X - \sum_{r=1}^d S^r \otimes P^r )^{-1} \ \eta \otimes v
\=
\eta \otimes (X-a)^{-1} v ,
\]
as one can verify by applying $
( I_{\h} \otimes X - \sum_{r=1}^d S^r \otimes P^r )$ to both sides.
So the right-hand side of (\ref{eqg2})
is 
\[
\la \eta, \xi \ra \la (X-a)^{-1} v , v \ra , \]
which is the same as the left-hand side of (\ref{eqg2}).
If $S$ has a spanning set of eigenvectors, our claim is proved.
If it does not, one can approximate it in norm by a $d$-tuple that does, and
the claim follows by continuity.
\bs
Now let $S$ and $T$ be $d$-tuples of commuting self-adjoint operators with 
$\sigma(S) \cup \sigma(T) \subset E$ and $\Delta := T - S \geq 0$. Let $$
R^r(t) \= (1-t) S^r + t T^r, \qquad 1 \leq r \leq d .$$
Then for $t$ in the range $(0,1)$, the $d$-tuple
$R(t)$ will consist of self-adjoint operators that need not commute
with each other.
Nonetheless,
letting
\[ Y(t) \=
\left( I_{\h} \otimes X - R(t) \odot I \right)
\]
then
$
R_v^*
Y(t)^{-1}
R_v $ 
still makes sense by Lemma~\ref{lemg1}.
Moreover,
\[
\frac{d}{dt} Y(t)^{-1} \=
Y(t)^{-1} (\Delta \odot I) \ Y(t)^{-1},
\]
and so is positive.
Therefore
\beq
F(T) - F(S) &\=& R_v^*Y(1)^{-1} R_v - R_v^*Y(0)^{-1} R_v \\
                        &=&
 R_v^* \int_0^1 \frac{d}{dt}Y(t)^{-1} dt R_v  \\
                        &=&
 \int_0^1 R_v^* Y(t)^{-1} (\Delta \odot I)
Y(t)^{-1} R_v dt \\
                        &\geq & 0.
\eeq
\ep
\bl
\label{lemg1}
Let $a^r < b^r, \ 1 \leq r \leq d$, and let 
$X$ be a densely defined self-adjoint operator on a graded Hilbert
space $\M \= \oplus_{r=1}^d \M^r $. Suppose that for every $t$ in $(0,1)$,
the point $\l_t = (1-t)a + t b$ is not in the $\mu$-spectrum of
$X$. Let $S = (S^1,\dots,S^d)$ be a $d$-tuple of bounded self-adjoint operators
on a Hilbert space $\h$,
with $\sigma(S^r) \subset (a^r,b^r)$ for each $r$.
Then $I \otimes X - \sum_{r=1}^d S^r \otimes P^r$ is invertible, with a bounded inverse.
\el
\bp
First, suppose $a^r = -1$ and $b^r = 1$
for each $r$. Then $(-1,1) \cap \sigma (X) $ is empty, and so 
is $(-1,1) \cap \sigma (I \otimes X)$.
So
$\| I_\h \otimes X \xi \| \geq \| \xi \| $ for every $\xi$ in $\h \otimes\M$.
But if $\sigma(S^r) \subset ( -1,1)$ for each $r$, the operator
$\sum S^r \otimes P^r$ has norm less than one. Therefore
$I \otimes X - \sum_{r=1}^d S^r \otimes P^r$ is invertible.

In the general case, let $m^r$ be the midpoint and $c^r$ half the length of the
interval $(a^r,b^r)$, so $a^r = m^r - c^r,\ b^r = m^r + c^r$.
Let 
\[
Y \= 
\left( \sum_{r=1}^d \frac{1}{\sqrt{c^r}} I \otimes P^r \right)
\left(
I \otimes X - \sum_{r=1}^d m^r I_\h \otimes P^r
\right)
\left( \sum_{r=1}^d \frac{1}{\sqrt{c^r}} I \otimes P^r \right) .
\]
If $\sigma_\mu$ denotes the $\mu$-spectrum,
\[
\sigma_\mu(Y) = c^{-1} (\sigma_\mu(X) - m),
\]
and hence the point $(1-t)(-\one) + t\one$ lies in the
$\mu$-resolvent set of $Y$ for $0<t<1$.  Let $T^r =
(1/c^r)(S^r-m^rI_\mathcal{H})$.  Then $T^r$ is a strict contraction, and
so, by the previous case, $Y - \sum T^r\otimes P^r$ is invertible.
As
\beq
\lefteqn{
Y - \sum_{r=1}^d T^r \otimes P^r}\\
&
\=
\left( \sum \frac{1}{\sqrt{c^r}} I \otimes P^r \right)
\left(
I \otimes X - \sum S^r \otimes P^r
\right)
\left( \sum \frac{1}{\sqrt{c^r}} I \otimes P^r \right) ,
\eeq
we get the desired result.
\ep

We can now prove a global result for rational functions of two variables.
\bt
\label{thmg2}
Let $F$ be a rational function of two variables. Let $\Gamma$ be the zero-set of the denominator
of $F$. Assume $F$ is real-valued on $\R^2 \setminus \Gamma$.
Let $E$ be an open rectangle in $\R^2 \setminus \Gamma$. Then $F$ is globally operator
monotone on $E$ if and only if $F$ is in $\L(E)$.
\et
\bp
Necessity follows from Theorem~\ref{thmloa}.

For sufficiency, by Lemma~\ref{lemg2} it is sufficient to prove the theorem
for $F_{t_n}$ with $t_n \searrow 0$, where $F_t = \rho_t \circ F \circ \rho_t$.  
Suppose the degree of $F$ is $n^1$ in $z^1$ and $n^2$ in $z^2$.
Let $\phi = \beta \circ F \circ \alpha$. 
By a result of G.~Knese \cite{kn08ub},
there is a model for $\phi$ in a Hilbert space $\M = \M^1 \oplus \M^2$ with
${\rm dim}(\M^r) = n^r$ for $r=1$ and $2$;  see also the paper
\cite{bsv05} by J.A.~Ball, C.~Sadosky and V.~Vinnikov.

Accordingly, in Theorem \ref{thmd2}, we obtain a realization of $F_t$ on
$\M^1\oplus\M^2$ of the form (\ref{eqd2}); since $\M$ is
finite-dimensional, the vector $v$ is in the domain of $X$.  By Theorem
\ref{thmd4} (v)$\Rightarrow$(ii), for some $v_1\in\M$
$$
F_t(z) \= C + \la (X-z)^{-1} v_1, v_1 \ra_{\M},
$$
where
${\rm dim}(\M^r) = n^r$ for $r=1$ and $2$,
and $F_t(\i, \i) = C < \i$. 
Then the pole-set $\Gamma_t$ of $F_t$ is contained in the zero-set
of ${\rm det}(X-z)$ which is a rational function of degree $({\rm dim}(\M^1),
{\rm dim}(\M^2))$. As these
two algebraic sets have the same degree, they must be equal.
So the $\mu$-resolvent of $X$ is $\R^2 \setminus \Gamma_t$, and now the result
follows from Theorem~\ref{thmg1}.
\ep
Let $\rho_t$ be as in (\ref{eqd15}). The following lemma is elementary.
\bl
\label{lemg2}
Let $t > 0$. Let $U$ be an open set in $\R^d$. Then:

(i) The function $F$ is globally operator monotone on $U \cap (-1/t,\i)^d$ if and only
if $F \circ \rho_t$ is globally operator monotone on $\rho_t^{-1} (U) \cap (-\i,1/t)^d$.

(ii) The function $F$ is globally operator monotone on $U \cap F^{-1}(-1/t,1/t) $
if and only if $\rho_t \circ F$ is globally operator monotone on the same set.
\el

What happens to Theorem~\ref{thmg2} in $d \geq 3$ variables?
It is still true that rational L\"owner functions have finite-dimensional
models \cite{colwer99, bsv05}. However, a recent example of
Knese \cite{kn10a} shows that the minimal dimension $n^r$ needed may be strictly greater
than the degree of $F$ in $z^r$.
So we cannot rule out the possibility that the $\mu$-spectrum of $X$ 
contains
some other algebraic sets in $\R^d$ than just the zero set of the denominator of $F$.  

We now give an example of a non-rational globally operator
monotone function.
\begin{example}
For $0 \leq s \leq 1/2$, the function $(z^1 z^2)^s$ is operator monotone
on $(0,\i) \times (0,\i)$.

Indeed, if $ (0,0) < (A^1, A^2) \leq (B^1, B^2)$ and $s$ is between $0$ and $1/2$,
then $$
\|  (A^r)^s (B^r)^{-s} \| \ \leq \ 1 \qquad {\rm for \ } r=1,2.
$$
Therefore the norm of 
\[
 (B^1)^{-s} (A^1)^s (A^2)^s (B^2)^{-s} 
\]
is less than or equal to $1$, so the largest eigenvalue  is less than or equal to $1$,
and therefore the largest eigenvalue of
\be
\label{eqz4}
(B^2)^{-s/2} (B^1)^{-s/2} (A^1)^s (A^2)^s (B^1)^{-s/2} (B^2)^{-s/2} 
\ee
is also less than or equal to $1$.
But (\ref{eqz4}) is self-adjoint, so less than or equal to the identity.
Therefore
\[
(A^1 A^2)^s \ \leq \ (B^1 B^2)^s .
\]

We do not know if $(z^1 z^2)^s$ can be approximated by rational functions in the
L\"owner class.

\end{example}

Let us close with some questions.

$\bullet$ 
Is Theorem~\ref{thmg2} true for rational 
functions of more than 2 variables?

$\bullet$ 
Can $E$ be an arbitrary open set in Theorem~\ref{thmg1}?

$\bullet$ 
Is every function in $\L(E)$ globally operator monotone on $E$?

$\bullet$ 
Is every function in $\lne$ $M_n$-monotone on $E$?